\documentclass{article}
\usepackage{amsmath}
\usepackage{amsthm}
\usepackage{amsfonts}
\usepackage{amscd}
\newcommand{\ZZ}{\mathbb{Z}}
\newcommand{\RR}{\mathbb{R}}
\newtheorem{theorem}{Theorem}
\newtheorem{lemma}[theorem]{Lemma}
\newtheorem{cor}[theorem]{Corollary}
\theoremstyle{definition}
\newtheorem{defn}[theorem]{Definition}
\newtheorem{example}[theorem]{Example}
\theoremstyle{remark}
\newtheorem*{acn}{Acnowledgements}
\numberwithin{equation}{section}
\numberwithin{theorem}{section}
\DeclareMathOperator{\Aut}{Aut}
\DeclareMathOperator{\Stab}{Stab}
\DeclareMathOperator{\per}{per}
\DeclareMathOperator{\imm}{imm}
\DeclareMathOperator{\Mat}{Mat}
\DeclareMathOperator{\sgn}{sgn}

\begin{document}

\title{Enumeration of walks on lattices. I}

\author{
Aleksandrs Mihailovs\\
Department of Mathematics\\
University of Pennsylvania\\
Philadelphia, PA 19104-6395\\
mihailov@math.upenn.edu\\
http://www.math.upenn.edu/$\sim$mihailov/
}

\date{\today}
\maketitle

\begin{abstract}
This work develops a methodical approach to counting of 
walks on cartesian products, biproducts, symmetric  and exterior 
powers  and bipowers, Schur operations, coverings  and semicoverings 
of weighted graphs. For weight  and root lattices of semisimple 
Lie algebras, this approach allows us to compute various combinatorial 
 and representation-theoretical constants, in particular, the number of 
plane symplectic wave graphs with given number of vertices.  
\end{abstract}    
\setlength{\baselineskip}{1.5\baselineskip}

\section{Introduction}

If it is not specified, a graph here means the graph in the most general sense, i.\ e.\ 
it can contain multiple edges, loops and some edges can be directed. The only restriction 
is that all graphs considered here, are locally finite, 
i.\ e.\ have a finite degree of each vertex. 
We denote $E(G)$  and $V(G)$ the sets of edges  and vertices of the graph $G$, as usual. 
Then 
\begin{equation}
E(G)=E_d(G)\amalg E_u(G),
\end{equation}
where $E_d(G)$  and $E_u(G)$ are the sets of the directed or undirected edges, correspondingly.
All loops are supposed to be directed.
Adjacency matrix $A(G)$ is the $V(G)\times V(G)$ matrix such that $A(G)_{uv}$ is the number 
of edges (both directed  and undirected) from $u$ to $v$.  

A graph morphism $\alpha: G_1\rightarrow G_2$ is a pair of maps 
\begin{equation}
\alpha_V: V(G_1)\rightarrow V(G_2),\quad \alpha_E: E(G_1)\rightarrow E(G_2)
\end{equation}
such that for any edge $e\in E(G_1)$ from $u$ to $v$ the edge $\alpha(e)\in E(G_2)$ is 
an edge from $\alpha(u)$ to $\alpha(v)$ of the same type, directed or undirected, as $e$. 
Graphs  and their morphisms form a category.

A weighted graph is a graph $G$ with the weight function $w: E(G)\rightarrow Z$, where 
$Z$ is a ring of weights. $Z$ can be any commutative ring and we suppose 
that this ring is fixed during this article. Weight matrix $W(G)$ of a weighted graph $G$ 
is the $V(G)\times V(G)$ $Z$-matrix such that 
\begin{equation}
W(G)_{uv}=\sum_{\substack{\text{$e$ from}\\ \text{$u$ to $v$}}}w(e) .
\end{equation}  

A weighted graph morphism $\alpha: G_1\rightarrow G_2$ is a graph morphism from $G_1$ to 
$G_2$ such that $w(\alpha(e))=w(e)$ for any edge $e\in E(G_1)$. Weighted graphs  and their 
morphisms form a category. 

If $Z$ is an ordered ring, i.\ e.\ its underlying set is ordered such that $a>0, b>0$ implies 
$a+b>0$  and $ab>0$ and if $w(e)>0$ for all edges $e\in E(G)$, then we'll say that $G$ is a 
positive weighted graph. Positive weighted graphs  and their morphisms form a category 
as well.

We can consider each graph $G$ as a positive weighted graph with $Z=\ZZ$  and $w\equiv 1$. 
In that case $A(G)=W(G)$.

We'll call a weighted graph {\em reduced} iff it is a graph without multiple edges, such 
that its weight function is nonzero for all its edges.  Reduced graphs  and their morphisms 
form a category, too. 

Cartesian product of (weighted) graphs $G_1\times G_2$ is a graph with 
\begin{equation}
V(G_1\times G_2)=V(G_1)\times V(G_2), \quad E(G_1\times G_2)=(V(G_1)\times E(G_2))\amalg 
(E(G_1)\times V(G_2))  
\end{equation}
supposing that for edge $e\in E(G_2)$ from $u$ to $v$, the edge $(v_1,e)$ is an edge from 
$(v_1,u)$ to $(v_1,v)$ of the same type, directed or undirected, as $e$ and $w(v_1,e)=
w(e)$ and analogously for the edges of another type. 

\begin{example}
Let $G=P_2$:

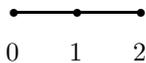
\begin{picture}(410, 25) 
\put(187,18){\circle*{3}}
\put(211,18){\circle*{3}}
\put(235,18){\circle*{3}}
\put(184,0){$0$}
\put(208,0){$1$}
\put(232,0){$2$}
\put(188,18){\line(1,0){22}}
\put(212,18){\line(1,0){22}}
\end{picture}

Then $P_2\times P_2$ looks as follows:

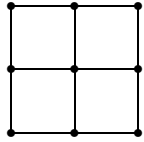
\begin{picture}(410, 52)
\put(187,50){\circle*{3}}
\put(211,50){\circle*{3}}
\put(235,50){\circle*{3}}
\put(187,26){\circle*{3}}
\put(211,26){\circle*{3}}
\put(235,26){\circle*{3}}
\put(187,2){\circle*{3}}
\put(211,2){\circle*{3}}
\put(235,2){\circle*{3}}

\put(188,50){\line(1,0){22}}
\put(212,50){\line(1,0){22}}
\put(188,26){\line(1,0){22}}
\put(212,26){\line(1,0){22}}
\put(188,2){\line(1,0){22}}
\put(212,2){\line(1,0){22}}

\put(187,3){\line(0,1){22}}
\put(211,3){\line(0,1){22}}
\put(235,3){\line(0,1){22}}
\put(187,27){\line(0,1){22}}
\put(211,27){\line(0,1){22}}
\put(235,27){\line(0,1){22}}
\end{picture}
\end{example}

Biproduct of (weighted) graphs $G_1\times_2 G_2$ is a graph with 
\begin{equation}
V(G_1\times_2 G_2)=V(G_1)\times V(G_2), \quad E(G_1\times_2 G_2)\simeq 
(E(G_1)\times E(G_2))\amalg (E_u(G_1)\times E_u(G_2))
\end{equation}
supposing that for edges $e_1\in G_1$ from $u_1$ to $v_1$  and $e_2\in G_2$ from 
$u_2$ to $v_2$, at least one of which is directed, the edge $(e_1,e_2)$ 
is a directed edge from $(u_1,v_1)$ to $(u_2,v_2)$, 
 and $w(e_1,e_2)=w(e_1)w(e_2)$; actually if one of those edges was undirected, considering it 
as an edge from $u$ to $v$, or from $v$ to $u$, we got two undirected edges in the biproduct;
for both undirected edges we got two undirected edges: 
from $(u_1,v_1)$ to $(u_2,v_2)$  and from $(u_1,v_2)$ to $(u_2,v_1)$, both of the same weight 
$w(e_1)w(e_2)$.

Below, in Example \ref{exbiprod} you could find a pictures of  
$P_2\times_2P_2$.

For unweighted either directed or undirected graphs biproduct 
$G=G_1\times_2G_2$ together with the canonical projections 
\begin{equation}\label{eqpr2}
\begin{split}
&\Pr\! {_1}: G\rightarrow G_1,\quad V(G)\ni(v_1,v_2)\mapsto v_1\in G_1,\\
&\Pr\! {_2}: G\rightarrow G_2,\quad V(G)\ni(v_1,v_2)\mapsto v_2\in G_2,
\end{split}
\end{equation}
projecting an edge from $(u_1,v_1)$ to $(u_2,v_2)$ onto the edges from $u_1$ to $u_2$  and 
from $v_1$ to $v_2$, correspondingly, is a product in the category-theoretical sense. 

\section{Cartesian products}\label{sec1}

\begin{defn}
Denote $c_n(G,u \rightarrow v)$ the count of walks of length $n$ starting in the vertex $u$ 
 and ending in the vertex $v$ of the graph $G$. Set 
\begin{equation}\label{1}
c_0(G,u \rightarrow v)=\delta_{uv}.
\end{equation}
Also denote 
\begin{equation}\label{2}
c(G,u \rightarrow v)=\sum_{n=0}^{\infty}\frac{c_n(G,u \rightarrow v)}{n!}t^n
\end{equation} 
the exponential generating function.
\end{defn} 

\begin{theorem}\label{thm1}
For a Cartesian product $G_1\times G_2$, 
\begin{equation}\label{3}
c(G_1\times G_2,(u_1,u_2)\rightarrow (v_1,v_2))=c(G_1,u_1\rightarrow v_1)
\thinspace c(G_2,u_2\rightarrow v_2) .
\end{equation}
\end{theorem}
\begin{proof}
By the definition of Cartesian product, 
for each edge $((a_1,a_2), (b_1,b_2))$ of 
$G_1\times G_2$, we have either $a_1=b_1$  and $(a_2,b_2)$ is an edge of $G_2$, or 
$(a_1,b_1)$ is an edge of $G_1$  and $a_2=b_2$. Choose a walk of length $n$ on the 
Cartesian product.
If $k$ of its edges are of the second kind, then the other $n-k$ edges 
are of the first kind. In this case $k$ edges of type $(a_1,b_1)$ form a walk 
of length $k$ on $G_1$ and $n-k$ edges of type $(a_2,b_2)$ form 
a walk of length $n-k$ on $G_2$. From the other side, for 
each pair of walks of length $k$ on $G_1$  and of length $n-k$ on $G_2$, we can 
construct a walk on $G_1\times G_2$ by $\binom{n}{k}$ ways, moving on 
$k$ arbitrarily chosen steps between $1$  and $n$ along  
the corresponding edge of the second kind and on the other $n-k$ steps---along 
the corresponding edge of the first kind. Hence  
\begin{equation}\label{4}
c_n(G_1\times G_2, (u_1,u_2)\rightarrow (v_1,v_2))=\sum_{j=0}^n\binom{n}{j}
c_j(G_1,u_1\rightarrow v_1)c_{n-j}(G_2,u_2\rightarrow v_2) .
\end{equation}
But it means the same as (\ref{3}).
\end{proof}

\begin{defn}
Denote $C(G)$ the $V(G)\times V(G)$ matrix with 
\begin{equation}\label{5}
C(G)_{uv}=c(G,u \rightarrow v) .
\end{equation}
\end{defn}

\begin{cor}\label{cor1}
For a Cartesian product $G_1\times G_2$,
\begin{equation}\label{6}
C(G_1\times G_2)=C(G_1)\otimes C(G_2) .
\end{equation}
\end{cor}
\begin{proof}
It immediately follows from Theorem \ref{thm1}  and the definition of a tensor product.
\end{proof}

\begin{lemma}\label{lem1}
$C(G)=e^{A(G)t}$.
\end{lemma}
\begin{proof}
By the definition of an exponent, 
\begin{equation}\label{7}
e^{A(G)t}=\sum_{n=0}^{\infty}\frac{(A(G))^n}{n!}t^n .
\end{equation}
So, we have to check that 
\begin{equation}\label{8}
((A(G))^n)_{uv}=c_n(G,u \rightarrow v) .
\end{equation}
But 
\begin{equation}\label{9}
((A(G))^n)_{uv}=\sum A(G)_{u_1v_1}\dots A(G)_{u_nv_n} ,
\end{equation}
where $(u_1,v_1), \dots, (u_n,v_n)$ is a walk between $u$  and $v$. Each summand in 
(\ref{9}) equals 1, that proves (\ref{8}).
\end{proof}

Formulas (\ref{8}--\ref{9}) suggest following 
\begin{defn}\label{def3}
For a weighted graph $G$, denote 
\begin{equation}\label{10} 
c_n(G,u \rightarrow v)=\sum w(u_1,v_1)\dots w(u_n,v_n)
\end{equation}
where $(u_1,v_1), \dots, (u_n,v_n)$ is a walk between $u$  and $v$ and
$w(u,v)$ is the weight of the edge $(u,v)$. Let $Z\{ t\} $ be $Z$-algebra of 
exponential series 
\begin{equation}
\sum_{n=0}^{\infty}z_n\frac{t^n}{n!}, \quad z_n\in Z .
\end{equation}
Denote $c(G,u \rightarrow v)\in Z\{ t\} $  and $C(G)$ by 
(\ref{2})  and (\ref{5}) as well. 
\end{defn}

\begin{lemma}\label{lem2}
For a weighted graph $G$, 
\begin{equation}\label{11}
C(G)=e^{W(G)t}.
\end{equation}
\end{lemma}
\begin{proof}
See the proof of Lemma \ref{lem1}.
\end{proof}

\begin{lemma}\label{lem3}
For a Cartesian product of graphs $G_1\times G_2$, 
\begin{equation}\label{12}
A(G_1\times G_2)=A(G_1)\otimes I_{V(G_2)}+I_{V(G_1)}\otimes A(G_2) ,
\end{equation}
where $I$ is an identity matrix.
\end{lemma}
\begin{proof}
$G_1\times G_2$ has two kinds of edges, see the proof of Theorem \ref{thm1}.
Each kind conforms to the corresponding summand in (\ref{12}).
\end{proof}

\begin{defn}\label{defW}
For matrices $W_1$  and $W_2$, denote 
\begin{equation}\label{14}
W_1\uplus W_2=W_1\otimes I+I\otimes W_2 ,
\end{equation}
where $I$ is an identity matrix of the size of $W_2$ first time and of the size of 
$W_1$ second time.
\end{defn}   

\begin{lemma}\label{lem4}
For a Cartesian product of weighted graphs $G_1\times G_2$,
\begin{equation}\label{15}
W(G_1\times G_2)=W(G_1)\uplus W(G_2) .
\end{equation}
\end{lemma}
\begin{proof}
The same as for Lemma \ref{lem3}. 
\end{proof}

\begin{theorem}\label{thm2}
For a Cartesian product of weighted graphs $G_1\times G_2$, 
\begin{equation}\label{16}
c(G_1\times G_2,(u_1,u_2)\rightarrow (v_1,v_2))=c(G_1,u_1\rightarrow v_1)
\thinspace c(G_2,u_2\rightarrow v_2) .
\end{equation}
\end{theorem}
\begin{proof}The same as for Theorem \ref{thm1}.
\end{proof}

\begin{cor}\label{cor2}
For a Cartesian product of weighted graphs $G_1\times G_2$,
\begin{equation}\label{17}
C(G_1\times G_2)=C(G_1)\otimes C(G_2) .
\end{equation}
\end{cor}
\begin{proof}
It immediately follows from Theorem \ref{thm2}  and the definition of a tensor product.
\end{proof}

\begin{cor}\label{cor3}
For any locally finite, i.\ e.\ with a finite number of nonzero elements in each row 
 and each column, matrices $W_1$  and $W_2$, 
\begin{equation}\label{18}
e^{(W_1\uplus W_2)t}=e^{W_1t}\otimes e^{W_2t} .
\end{equation}
\end{cor}
\begin{proof}
It follows from Lemmas \ref{lem2}  and \ref{lem4}  and Corollary \ref{cor2}.
\end{proof}

\section{Biproducts}

\begin{defn}
Denote $\overrightarrow{G}$ the directed (weighted) graph obtained from a 
(weighted) graph $G$ by preserving all its 
vertices  and directed edges (with their weights), 
 and by the replacement of each of its undirected edges 
by a pair of opposite directed edges (of the same weight) between the same vertices.
\end{defn}

\begin{lemma}\label{lembi1}
$A(\overrightarrow{G})=A(G)$  and $W(\overrightarrow{G})=W(G)$.
\end{lemma}
\begin{proof}
It follows immediately from the definitions of the matrices $A(G)$  and $W(G)$.
\end{proof}

\begin{cor}\label{corbi1}
$C(\overrightarrow{G})=C(G)$.
\end{cor}
\begin{proof}
It follows from Lemmas \ref{lem1}, \ref{lem2} and \ref{lembi1}.
\end{proof}

\begin{lemma}\label{lembi2}
$\overrightarrow{G_1\times_2G_2}\simeq\overrightarrow{G_1}\times_2\overrightarrow{G_2}$.
\end{lemma} 
\begin{proof}
Again it follows immediately from the definitions.
\end{proof}

\begin{defn}\label{dcirc}
Denote $\circ$ the coefficient-wise product in $Z\{ t\} $:
\begin{equation}
\left( \sum_{n=0}^{\infty}x_n\frac{t^n}{n!}\right) \circ 
\left( \sum_{n=0}^{\infty}y_n\frac{t^n}{n!}\right) = 
\sum_{n=0}^{\infty}x_ny_n\frac{t^n}{n!} .
\end{equation}
\end{defn}

\begin{theorem}\label{thbi1}
For a biproduct $G_1\times_2 G_2$ of (weighted) graphs, 
\begin{gather}\label{eqthbi11}
c(G_1\times_2 G_2, (u_1,u_2)\rightarrow (v_1,v_2))=c(G_1, u_1\rightarrow v_1)
\circ c(G_2, u_2\rightarrow v_2) ,\\
c_n(G_1\times_2 G_2, (u_1,u_2)\rightarrow (v_1,v_2))=c_n(G_1, u_1\rightarrow v_1)
\thinspace c_n(G_2, u_2\rightarrow v_2) .
\label{eqthbi12}
\end{gather}
\end{theorem}
\begin{proof}
By Corollary \ref{corbi1}  and Lemma \ref{lembi2}, we can suppose that the graphs 
$G_1$  and $G_2$ are directed. Consider the canonical projections $\Pr\!{_1}$  and 
$\Pr\!{_2}$ defined in \eqref{eqpr2}. They project each walk on $G_1\times_2G_2$ 
onto two walks of the same length on $G_1$  and $G_2$. Conversely, each pair of walks 
of equal lengths on $G_1$  and $G_2$ produce a walk on $G_1\times_2G_2$ of the same length, 
where we choose the product of the first edges on the first step, the second edges on 
the second step and so on. Thus, we have a bijection between the corresponding sets 
of walks on $G_1\times_2G_2$  and pairs of walks of the same length on $G_1$  and $G_2$.
That proves \eqref{eqthbi12} and \eqref{eqthbi11} follows from \eqref{eqthbi12} for 
all $n$. The weight of a walk on $G_1\times_2G_2$ equals the product of the weights 
of the corresponding walks on $G_1$  and $G_2$ as well. 
\end{proof}

\begin{defn}
We'll call a graph $G$ bipartite if the set of its vertices is divided on two 
nonintersecting parts
\begin{equation}
V(G)=V_0(G)\amalg V_1(G)
\end{equation}
such that there are no edges between vertices lying in the same part, i.\ e.\ each 
edge of $G$ has an end in $V_0(G)$  and another one---in $V_1(G)$. In particular, bipartite 
graphs can't contain loops.
\end{defn}

\begin{defn}
Let $G_1$  and $G_2$ be bipartite graphs. We'll say that a vertice $(v_1,v_2)$ of their 
biproduct $G_1\times_2G_2$, is {\em even}, if either 
\begin{equation}
v_1\in V_0(G_1),\thickspace v_2\in V_0(G_2),\quad \text{or}\quad  
v_1\in V_1(G_1),\thickspace v_2\in V_1(G_2),
\end{equation}
 and {\em odd} otherwise, i.\ e.\ if either 
\begin{equation}
v_1\in V_0(G_1),\thickspace v_2\in V_1(G_2),\quad \text{or}\quad  
v_1\in V_1(G_1),\thickspace v_2\in V_0(G_2).
\end{equation}
\end{defn}

\begin{lemma}\label{lembi3}
Each edge of a biproduct of bipartite graphs connects the 
vertices of the same parity, i.\ e.\ its ends are either both even, or both odd.
\end{lemma}
\begin{proof}
For an edge connecting the vertices $(u_1,u_2)$  and $(v_1,v_2)$ of $G_1\times_2G_2$,
if 
\begin{equation}
u_1\in V_a(G_1),\quad u_2\in V_b(G_2),
\end{equation}
then
\begin{equation}
v_1\in V_{1-a}(G_1),\quad v_2\in V_{1-b}(G_2),
\end{equation}
 and
\begin{equation}
a+b\equiv (1-a)+(1-b) \mod 2,
\end{equation}
that means that the vertices $(u_1,u_2)$  and $(v_1,v_2)$ are both even or odd 
simultaneously.
\end{proof}

\begin{defn}
For a bipartite (weighted) graphs $G_1$  and $G_2$ denote $G_1\times_0G_2$  and $G_1\times_1G_2$ 
the complete (weighted) subgraphs of the biproduct $G_1\times_2G_2$, spanned by the even or odd 
vertices, correspondingly. We'll call these graphs an even  and an odd product of $G_1$  and 
$G_2$, correspondingly.
\end{defn}

\begin{example}\label{exbiprod}
Let $G=P_2$ with $V_0(P_2)=\{ 0, 2\} $, $V_1(P_2)=\{ 1\} $:

\begin{picture}(410, 25) 
\put(187,18){\circle*{3}}
\put(211,18){\circle*{3}}
\put(235,18){\circle*{3}}
\put(184,0){$0$}
\put(208,0){$1$}
\put(232,0){$2$}
\put(188,18){\line(1,0){22}}
\put(212,18){\line(1,0){22}}
\end{picture}

Then $G\times_0G$  and $G\times_1G$ look as follows:

\begin{picture}(410, 75)
\put(92,66){\circle*{3}}
\put(140,66){\circle*{3}}
\put(296,66){\circle*{3}}
\put(116,42){\circle*{3}}
\put(272,42){\circle*{3}}
\put(320,42){\circle*{3}}
\put(92,18){\circle*{3}}
\put(140,18){\circle*{3}}
\put(296,18){\circle*{3}}

\put(93,65){\line(1,-1){22}}
\put(117,43){\line(1,1){22}}
\put(93,19){\line(1,1){22}}
\put(117,41){\line(1,-1){22}}
\put(273,43){\line(1,1){22}}
\put(297,65){\line(1,-1){22}}
\put(273,41){\line(1,-1){22}}
\put(297,19){\line(1,1){22}}

\put(97,0){$P_2\times_0P_2$}
\put(277,0){$P_1\times_1P_2$}
\end{picture}
\end{example}

\begin{cor}\label{corbi2}
For a bipartite (weighted) graphs $G_1$  and $G_2$, 
\begin{gather}\label{eqcorbi21}
c(G_1\times_0 G_2, (u_1,u_2)\rightarrow (v_1,v_2))=c(G_1, u_1\rightarrow v_1)
\circ c(G_2, u_2\rightarrow v_2) ,\\
c(G_1\times_1 G_2, (u_1,u_2)\rightarrow (v_1,v_2))=c(G_1, u_1\rightarrow v_1)
\circ c(G_2, u_2\rightarrow v_2) .
\label{eqcorbi22}
\end{gather}
\end{cor}
\begin{proof}
It is a corollary of Theorem \ref{thbi1}  and Lemma \ref{lembi3}.
\end{proof}

\section{Homogeneous graphs}\label{sec1a}

\begin{defn}\label{TI}
A pair of (weighted) graph $G$  and a group of its automorphisms $\Gamma\subseteq \Aut G$
forms a homogeneous (weighted) graph if for each two vertices $u, v$ of $G$
there is an automorphism $g\in\Gamma$ such that $g(u)=v$. If $\Gamma$ is an Abelian group,
we'll say that $G$ is an Abelian (weighted) graph and $\Gamma$ is its translation group.
\end{defn}

\begin{lemma}\label{lemhom}
Let $(u,u_1),\dots,(u,u_d)$ be the list of all the edges starting from a vertex 
$u$ of a (weighted) homogeneous graph $(G,\Gamma)$ ( and $w_1,\dots,w_d$ be their weights).
Then for any vertex $v$ of $G$ and for any $g\in\Gamma$ such that $v=g(u)$, the list 
of all the edges starting from $v$ is $(v,g(u_1)),\dots,(v,g(u_d))$ ( and they have the 
same weights $w_1,\dots,w_d$).
\end{lemma}
\begin{proof}
Because $g$ is an automorphism.
\end{proof}

\begin{defn}\label{defhom}
Let $\Gamma$ be a group, $H$ its subgroup and $g_1,\dots,g_d$ elements of 
$\Gamma$. Denote $(\Gamma, H$, $(w_1)g_1, \dots, (w_d)g_d)$ a (weighted) graph
with the set of vertices 
$\Gamma/H$, edges of which, starting in $g H$ are $(gH, gg_1H), \dots, 
(g H, gg_dH)$ (with weights $w_1, \dots, w_d$, correspondingly) for any $g H \in \Gamma/H$.
\end{defn} 

\begin{theorem}\label{lemhom1}
The pair $((\Gamma, H, (w_1)g_1, \dots, (w_d)g_d), \Gamma/H_{\Gamma})$ where 
\begin{equation}\label{etag}
H_{\Gamma}=\bigcap_{g\in \Gamma}gHg^{-1} 
\end{equation}
is a homogeneous (weighted)
graph for any admissible data.
Each (weighted) homogeneous graph $(G,\Gamma)$, edges of which (of weights $w_1,\dots,w_d$),
starting in its vertex $u$, 
are $(u,u_1),\dots,(u,u_d)$, is isomorphic to a graph 
$(\Gamma, H, (w_1)g_1, \dots, (w_d)g_d)$ with $H_{\Gamma}=0$
where $H\subseteq\Gamma$ 
is the stabilizer of a vertex $u$ and $u_1=g_1(u),\dots,u_d=g_d(u)$.
\end{theorem}
\begin{proof}
Elements of $ H_{\Gamma}$ act trivially on $\Gamma/ H$ and if 
$xg H=g H$, then $x\in g H g^{-1}$. It means that $x\in H_{\Gamma}$ iff $x$ acts 
trivially on $\Gamma/ H$. By definition \eqref{etag}, $H_{\Gamma}$ is a normal subgroup of 
$\Gamma$, so $\Gamma/ H_{\Gamma}\subseteq \Aut (\Gamma, H,\dots)$ 
 and it acts transitively on $\Gamma/ H$, because $\Gamma$ acts transitively. 
Everything remaining except the triviality of $ H_{\Gamma}$ follows from Lemma \ref{lemhom}. 
$ H_{\Gamma}\subseteq\Gamma\subseteq\Aut G$ must be trivial, because only the neutral 
element of $\Aut G$ acts trivially on $G$.
\end{proof}

\begin{theorem}\label{thcon}
Let $(G,\Gamma)$ be a connected homogeneous (weighted) graph and  $(u, g_1(u)), \dots, 
(u,$ $g_d(u))$ be
a list of all its edges, starting in a vertex $u$. Then $(G,\Gamma_0)$ is a homogeneous 
(weighted) graph as well, where $\Gamma_0$ is the subgroup of $\Gamma$, generated by 
$g_1, \dots, g_d$.
\end{theorem}
\begin{proof}
By Lemma 
\ref{lemhom}, each edge, starting from $\Gamma_0u$ has the form $(g(u), gg_i(u))$ with 
$g\in\Gamma_0$, hence it ends in $\Gamma_0u$ as well, because $gg_i\in\Gamma_0$. 
Analogously, the edges ending in $\Gamma_0u$, has the form $(g(u), gg_i(u))$ with 
$gg_i\in\Gamma_0$ and they start in $\Gamma_0u$ as well, because $g=(gg_i)g_i^{-1}\in
\Gamma_0$. It means that $\Gamma_0u$ contains the connected component of $u$, i.\ e.\  
\begin{equation}\label{G0u}
\Gamma_0u=G.
\end{equation}
\end{proof}

\begin{theorem}\label{thAb}
Each Abelian (weighted) graph with the translation group $\Gamma$ is isomorphic to a graph 
$(\Gamma, 0, (w_1)g_1,\dots,(w_d)g_d)$. If $G$ is connected, then elements  
$g_1,\dots,g_d$ generate $\Gamma$. 
\end{theorem}
\begin{proof}
We already know from Theorem \ref{lemhom1} that a graph is isomorphic to a
$(\Gamma, H,(w_1)g_1,\dots,(w_d)g_d)$ with $ H_{\Gamma}=0$. For an Abelian $\Gamma$ 
we have $ H_{\Gamma}= H$, hence $ H=0$. It means that for each vertex $v$ of $G$, there is 
a unique element $g\in G$ such that $v=gu$. Now, from \eqref{G0u}, $\Gamma=\Gamma_0$.
\end{proof} 

\begin{cor}\label{corAb}
The translation group of a connected Abelian (weighted) graph is finite generated, i.\ e.\ 
isomorphic to $\ZZ^r\oplus(\ZZ/q_1\ZZ)\oplus\dots\oplus(\ZZ/q_k\ZZ)$ where $q_1,\dots,q_k$ 
are some powers of prime numbers.
\end{cor}
\begin{proof}
It follows immediately from Theorem \ref{thAb}.
\end{proof}

\begin{lemma}\label{lemh}
For any vertices $u,v$ of a (weighted) homogeneous graph $G$  
 and an automorphism $g$ of $G$, 
\begin{gather}
c_n(G,u\rightarrow v)=
c_n(G,g(u)\rightarrow g(v))\\ 
c(G,u\rightarrow v)=
c(G,g(u)\rightarrow g(v)) . 
\end{gather}
\end{lemma}
\begin{proof} 
Applying $g$ to both vertices of each edge of a walk from $u$ to $v$, we get a walk from 
$g(u)$ to $g(v)$. 
Applying $g^{-1}$ to both vertices of each edge of a walk from $g(u)$ to $g(v)$, 
we get a walk from $u$ to $v$. Since these operations are inverse, we get a bijection 
between the corresponding sets.
\end{proof}

\begin{defn}\label{H}
For an automorphisms $g\in\Gamma$ of a (weighted) homogeneous graph 
$G=(\Gamma,H,\dots)$ denote 
\begin{gather}\label{18a}
c_n(G,g)=c_n(G,H\rightarrow gH)\\  
c(G,g)=c(G,H\rightarrow gH).
\label{18aa}\end{gather}
\end{defn} 

\begin{lemma}\label{lemH}
For any vertices $g_0H,g_1H$ of a (weighted) homogeneous graph 
$G=(\Gamma,H,\dots)$,
\begin{gather}
c_n(G,g_0H\rightarrow g_1H)=c_n(G,g_0^{-1}g_1) ,\\
c(G,g_0H\rightarrow g_1H)=c(G,g_0^{-1}g_1) .
\end{gather}
\end{lemma}
\begin{proof} 
It follows from Lemma \ref{lemh} with $u=H$, $v=g_0^{-1}g_1H$  and $g=g_0$. 
\end{proof}

\begin{theorem}\label{thtran}
Cartesian product of homogeneous (weighted) graphs is a homogeneous (weighted) graph, 
 and for graphs defined in Definition \ref{defhom},
\begin{equation}
c(G_1\times G_2, (g_1,g_2))=c(G_1, g_1)\thinspace c(G_2, g_2) .
\end{equation}
\end{theorem}
\begin{proof}
Cartesian product of automorphisps is an automorphism and we can use Theorems \ref{thm1} 
 and \ref{thm2}.
\end{proof} 

\begin{defn}
Let $Z$ be the ring of weights and $\Gamma$ be a group. Denote $Z(\Gamma)\{ t\}$
the $Z\{ t\}$-algebra of exponential series 
\begin{equation}\label{zg}
\sum_{n=0}^{\infty}z_n\frac{t^n}{n!}, \quad z_n\in Z(\Gamma), 
\end{equation}
with coefficients in the group algebra 
\begin{equation}\label{zg1} 
Z(\Gamma)=Z\otimes \ZZ(\Gamma) .
\end{equation}
\end{defn}

\begin{theorem}\label{thexp}
For $G=(\Gamma, 0, (w_1)g_1, \dots, (w_d)g_d)$, 
the matrix $C(G)$ is the matrix of right multiplication on  
$e^{((w_1)g_1+\dots+(w_d)g_d)t}$ in the topological basis $(g)_{g\in\Gamma}$ of 
the $Z\{ t\}$-algebra $Z(\Gamma)\{ t\}$.
\end{theorem}
\begin{proof}
By Lemma \ref{lemhom}, the adjacency (weight) matrix $A(G)$ (or $W(G)$ correspondingly) 
is the matrix $M$ of right 
multiplication on $((w_1)g_1+\dots+(w_d)g_d)$ in the given basis, i.\ e.\ for $g\in\Gamma$,
\begin{equation}\label{eqexp}
gM=\sum_{h\in\Gamma} M_{gh}h .
\end{equation}
Now, Theorem \ref{thexp} follows from Lemmas \ref{lem1}  and \ref{lem2}.
\end{proof}

\begin{cor}\label{corexp}
For $G=(\Gamma, 0, (w_1)g_1, \dots, (w_d)g_d)$, 
\begin{gather}\label{eqcorexp1}
((w_1)g_1+\dots+(w_d)g_d)^n=\sum_{g\in\Gamma}c_n(G,g)g ,\\
\label{eqcorexp2}
e^{((w_1)g_1+\dots+(w_d)g_d)t}=\sum_{g\in\Gamma}c(G,g)g .
\end{gather}
\end{cor}
\begin{proof}
\eqref{eqcorexp2} follows from Theorem \ref{thexp}  and \eqref{eqexp} 
for $M=C(G)$  and neutral $g$ in \eqref{eqexp}. Considering coefficients at 
$t^n$, we get \eqref{eqcorexp1}.
\end{proof}

\section{Composite walks}

\begin{defn}\label{through}Let $k\geq 2$. 
For vertices $v_0, v_1, \dots, v_k$ of a (weighted) graph $G$ denote 
$c_n(G,v_0 \rightarrow v_1 \rightarrow \dots \rightarrow v_k)$ 
the number of walks (the sum (\ref{10}) along walks) from $v_0$  coming 
successively through $v_1, \dots, v_{k-1}$ to $v_k$. 
Denote 
\begin{equation}\label{eqth}
c(G,v_0 \rightarrow v_1 \rightarrow \dots \rightarrow v_k)=\sum_{n=0}^{\infty}
\frac{c_n(G,v_0 \rightarrow v_1 \rightarrow \dots \rightarrow v_k)}{n!}t^n .
\end{equation}
\end{defn}

\begin{lemma}\label{lemh2}
For vertices $v_0,v_1,\dots,v_k$ of a (weighted) homogeneous graph $G$  
 and an automorphism $g$ of $G$, 
\begin{gather}
c_n(G,v_0\rightarrow v_1\rightarrow\dots\rightarrow v_k)=
c_n(G,g(v_0)\rightarrow g(v_1)\rightarrow\dots\rightarrow g(v_k)) ,\\ 
c(G,v_0\rightarrow v_1\rightarrow\dots\rightarrow v_k)=
c(G,g(v_0)\rightarrow g(v_1)\rightarrow\dots\rightarrow g(v_k)) . 
\end{gather}
\end{lemma}
\begin{proof} 
Applying $g$ to both vertices of each edge of a walk from $v_0$ coming successively through 
$v_1, \dots, v_{k-1}$ to $v_k$, we get a walk from 
$g(v_0)$ coming succesively through $g(v_1), \dots, g(v_{k-1})$ to $g(v_k)$. 
Applying $g^{-1}$ to both vertices of each edge of a walk from $g(v_0)$ to $g(v_k)$, 
we get a walk from $v_0$ to $v_k$. Since these operations are inverse, we get a bijection 
between the corresponding sets.
\end{proof}
\begin{defn}\label{H2}
For elements $g_1,\dots,g_k\in\Gamma$ of a (weighted) homogeneous graph 
$G=(\Gamma,H,\dots)$ denote 
\begin{gather}\label{18a2}
c_n(G,g_1,\dots,g_k)=c_n(G,H\rightarrow g_1H\rightarrow g_1g_2H\rightarrow\dots
\rightarrow g_1\dots g_kH) ,\\  
c(G,g_1,\dots,g_k)=c(G,H\rightarrow g_1H\rightarrow g_1g_2H\rightarrow\dots
\rightarrow g_1\dots g_kH) .
\label{18aa2}\end{gather}
\end{defn} 

\begin{lemma}\label{lemH2}
For vertices $g_0H,g_1H,\dots,g_kH$ of a (weighted) homogeneous graph 
$G=(\Gamma,H,\dots)$
\begin{gather}
c_n(G,g_0H\rightarrow g_1H\rightarrow\dots\rightarrow g_kH)=c_n(G,g_0^{-1}g_1,g_1^{-1}g_2,
\dots,g_{k-1}^{-1}g_k) ,\\
c(G,g_0H\rightarrow g_1H\rightarrow\dots\rightarrow g_kH)=c(G,g_0^{-1}g_1,g_1^{-1}g_2,
\dots,g_{k-1}^{-1}g_k) .
\end{gather}
\end{lemma}
\begin{proof} 
It follows from Lemma \ref{lemh2} with $v_i=g_0^{-1}g_iH$  and $g=g_0$. 
\end{proof}

\begin{theorem}\label{thfun}
Let $G=(\Gamma,H,\dots)$  and $F$ be a fundamental domain of the actions of some subgroup of
$\Stab H=$ $\{ \alpha\in\Aut G: \alpha(H)=H\}$ on $V(G)=\Gamma/H$. Then 
\begin{equation}
c(G,g_1,\dots,g_k)=c(G,|g_1|,\dots,|g_k|) ,
\end{equation}
where $|g|$ denotes an element of $\Gamma$ such that $|g|H\in F$  and $|g|H=\alpha_g(gH)$ 
for an automorphism $\alpha_g$ from the given subgroup of $\Stab H$.
\end{theorem}
\begin{proof}
For each $g\in\Gamma$, fix $\alpha_g$ satisfying the condition above. 
For each walk from $H$ coming successively through $g_1H,\dots,g_1\dots g_{k-1}H$ to  
$g_1\dots g_kH$ we can construct a walk from $H$ coming successively through 
$|g_1|, \dots, |g_1|\dots |g_{k-1}|H$ to $|g_1|\dots |g_k|H$. 
Walking from $H$ to $|g_1|H$, first we go along
the image under $\alpha_{g_1}$ of the initial 
part of the given walk to the first meeting of $g_1H$. Note that we'll meet $|g_1|H$ only 
at the end of this part of a walk, because $\alpha_{g_1}$ is an automorphism. Then 
we successively transform the part of the given walk from $g_1\dots g_{i-1}$ to the first 
(new) meeting of $g_1\dots g_iH$ into its image under $|g_1|\dots |g_{i-1}|\alpha_{g_i}$.
Again we'll meet $|g_1|\dots |g_i|$ only at the end, because we use an automorphism.
Conversely, for each walk from $H$ to $|g_1|\dots |g_k|H$ coming succesively through 
$|g_1|,\dots,|g_1|\dots|g_{k-1}|H$ we can construct a walk of the given type by inverse 
automorphisms of the corresponding parts. So, we have a bijection between the considered 
sets of walks.
\end{proof}

\begin{theorem}\label{thfun1}
If in the conditions of Theorem \ref{thfun} for any $f_1,f_2\in F$ each walk from $H$ to 
$f_1f_2H$ comes through $f_1H$, then 
\begin{equation}\label{eqfun1}
c(G,g_1,\dots,g_k)=c(G,|g_1|\dots|g_k|)
\end{equation}
\end{theorem}
\begin{proof}
For any $f\in F$, each walk from $fH$ to $ff_1f_2H$ comes through 
$ff_1H$, because multiplication on $f$ is an automorphism of $G$. It means that each walk 
from $H$ to $|g_1|\dots |g_k|H$ comes succesively through $g_1H, \dots, g_1\dots g_{k-1}H$. 
Therefore
\begin{equation}
c(G,|g_1|,\dots,|g_k|)=c(G,|g_1|\dots |g_k|) .
\end{equation}
Using Theorem \ref{thfun}, we get \eqref{eqfun1}.
\end{proof}

\section{Coverings  and semicoverings}

\begin{defn}\label{defcov}
We'll say that a (weighted) graph morphism $\pi: G_1\rightarrow G_2$ is a covering, 
iff it is a surjective local isomorphism, i.\ e.\ for each vertex $v\in G_2$ there is 
a vertex $u\in G_1$ such that 
\begin{equation}\label{pi}
v=\pi(u) ,
\end{equation}
 and each vertex $u\in G_1$ satisfying \eqref{pi} has the same count ( and weights)
of both incoming 
 and outcoming edges starting or ending in $u$ as the edges starting  and ending in $v$.
\end{defn}

\begin{theorem}\label{thcov}
Let $\pi: G_1\rightarrow G_2$ be a covering and $\pi(u_1)=u$, $\pi(v_1)=v$ for 
some vertices of $G_1$  and $G_2$. Then
\begin{equation}\label{thpi}
c(G_2,u\rightarrow v)=\sum_{v_j\in\pi^{-1}(v)} c(G_1,u_1\rightarrow v_j)=
\sum_{u_i\in\pi^{-1}(u)}c(G_1,u_i\rightarrow v_1).
\end{equation}
\end{theorem}
\begin{proof}
$\pi$ can't paste together edges starting or ending in the same vertex of $G_1$, because it 
is a local isomorphism. It means that distinct walks on $G_1$ projects under $\pi$ in distinct 
walks on $G_2$. Conversely, for each walk from $u$ to $v$ on $G_2$, choosing for the first  
edge some edge of $G_1$ starting in $u_1$  and projecting onto the given edge and then 
step by step choosing for each sequential edge of the walk on $G_2$ some edge of $G_1$ 
starting in the end of the edge of $G_1$ constructed on the previous step and projecting in 
the considered edge, we get a walk on $G_1$ from $u_1$ to some vertex $v_i$ projecting onto 
$v$. Therefore, we have a bijection between the set of walks on $G_2$ from $u$ to $v$  and the 
set of walks on $G_1$ from $u_1$ to vertices from $\pi^{-1}(v)$ and the first equality in 
\eqref{thpi} is proven. The equality between the first  and the last items of \eqref{thpi} 
can be proven analogously, or it follows from the proven first equality by changing
the directions of all edges.
\end{proof}

\begin{defn}\label{defmod}
For a graph $G$, denote $|G|$ the undirected graph with the same vertices as $G$, obtained from 
$G$ by saving all its undirected edges and a replacement of all directed edges by undirected 
ones having the same vertices. Also, for a graph morphism $\alpha: G_1\rightarrow G_2$ we 
denote $|\alpha|:|G_1|\rightarrow |G_2|$ a morphism that is the same as $\alpha$ on vertices 
 and undirected edges of $G_1$ and transfers an edge of $|G_1|$ obtained by a replacement of 
a directed edge $e$ of $G_1$ to the edge of $|G_2|$ obtained by a replacement of the directed 
edge $\alpha(e)$. 
\end{defn}

\begin{lemma}\label{lemmod}
$|.|$ is a functor from the category of graphs to the category of 
undirected graphs.
\end{lemma}
\begin{proof}
It is clear.
\end{proof}

\begin{lemma}\label{lemmod1}
The morphism $|\alpha|$ is a covering iff $\alpha$ is a covering.
\end{lemma}
\begin{proof}
One can easily check that both conditions, of a surjectivity  and of being a local 
isomorphism, are true or not simultaneously for $\alpha$  and $|\alpha|$.
\end{proof}  

\begin{defn}
We'll call a graph $G$ connected iff $|G|$ is connected. 
\end{defn}
 
\begin{theorem}\label{thcov1}
Let $G_1$ be a connected graph  and $\pi: G_1\rightarrow G_2$ be a covering. 
Then $G_2$ is connected as well and either the sets 
$\pi^{-1}(v)$ are infinite for each vertex $v\in G_2$, or all of them are finite 
 and contain an equal number of elements.
\end{theorem}
\begin{proof}
Morphism $\pi$ is surjective and
the image of a morphism of a connected graph is connected, so $G_2$ is connected.
By Lemma \ref{lemmod1}, we can suppose that the graphs $G_1$  and $G_2$ are undirected, 
or replace them by $|G_1|$  and $|G_2|$ otherwise. 
Now, if the set $\pi^{-1}u$ is finite  and contains $m_u>0$ elements for some $u\in G_2$, then 
for any $v\in G_2$ we obtain from Theorem \ref{thcov} that 
\begin{equation}\label{eqcov2}
\sum_{u_i\in\pi^{-1}(u)}\sum_{v_j\in\pi^{-1}(v)}c(G_1, u_i\rightarrow v_j)=
m_u\thinspace c(G_2, u\rightarrow v). 
\end{equation}  
All coefficients of the series in \eqref{eqcov2} are nonnegative if we consider 
regular walks, not weighted,
therefore we can change the order of summation:
\begin{equation}\label{eqcov3}
\sum_{u_i\in\pi^{-1}(u)}\sum_{v_j\in\pi^{-1}(v)}c(G_1, u_i\rightarrow v_j)=
\sum_{v_j\in\pi^{-1}(v)}\sum_{u_i\in\pi^{-1}(u)}c(G_1, u_i\rightarrow v_j)=
m_v\thinspace c(G_2, u\rightarrow v),
\end{equation}
where $m_v$ denotes the number of elements of $\pi^{-1}(v)$, that must be finite to ensure 
the converges of the series, because 
\begin{equation}\label{eqcov4}
c(G_2, u\rightarrow v)\neq 0
\end{equation}
since $G_2$ is connected. Compare \eqref{eqcov2}  and \eqref{eqcov3}, we get 
\begin{equation}\label{eqcov5}
m_u=m_v.
\end{equation}
\end{proof}

\begin{cor}\label{corcov}
Let $\pi: G_1\rightarrow G_2$ be a covering. Then for each connected component $G$ of $G_1$, 
the subgraph $\pi(G)\subseteq G_2$ is a connected component of $G_2$  and either the sets 
$\pi^{-1}(v)\cap G$ are infinite for each vertex $v\in \pi(G)$, or all of them are finite 
 and contain an equal number of elements.
\end{cor}
\begin{proof}
For each edge of $G_2$ starting or ending in a vertex $v\in\pi(G)$, we can find an edge of 
$G$ projecting onto it, because $\pi$ is a local isomorphism. It means that $\pi(G)$ is 
a connected component of $G_2$. The other follows from Theorem \ref{thcov1}.
\end{proof}

\begin{defn}
We'll say that a covering $\pi:G_1\rightarrow G_2$ is infinite-sheeted if for each vertex 
$v\in G_2$ the set $\pi^{-1}(v)$ is infinite, or $j$-sheeted if for each vertex $v\in G_2$ 
the set $\pi^{-1}(v)$ contains exactly $j$ elements.
\end{defn}

\begin{lemma}\label{lemsh}
Let $\Gamma$ be a group and $H\subseteq H_1\subseteq G$ be their subgroups. Then the 
morphism 
\begin{equation}\label{eqsh}
\begin{split}
\pi:(\Gamma, H, (w_1)g_1, \dots, (w_d)g_d)&\rightarrow (\Gamma, H_1, (w_1)g_1,\dots, 
(w_d)g_d),\\ gH\mapsto gH_1,\quad (gH,gg_kH)&\mapsto (gH_1, gg_kH_1)
\end{split} 
\end{equation}
is an $(H_1:H)$-sheeted covering.
\end{lemma}
\begin{proof}
First we have to check that \eqref{eqsh} defines a (weighted) graph morphism.
In fact, $gH=fH$ iff $f=gh$ for some $h\in H\subseteq H_1$, hence $gH_1=fH_1$ as well.
Thus, $\pi$ is a (weighted) graph morphism and it is a local isomorphism  and surjective 
by its construction, so it is a covering. Cosets projected onto $gH_1$ are $ghH$ with 
$h\in H_1$. There are exactly $(H_1:H)$ cosets of that type. 
\end{proof}

\begin{cor}
Let $G_1=(\Gamma, H, (w_1)g_1, \dots, (w_d)g_d)$, $G_2=(\Gamma, H_1, (w_1)g_1, \dots, 
(w_d)g_d)$, with $H\subseteq H_1$  and $\pi: G_1\rightarrow G_2$ be a covering, defined 
in \eqref{eqsh}. Then 
\begin{equation}
c(G_2, g)=\sum_{hH\in H_1/H}c(G_1,gh) ,
\end{equation}
 and, in particular, if $H=\{ 1\}$, then 
\begin{equation}
c(G_2,g)=\sum_{h\in H_1}c(G_1,gh).
\end{equation}
\end{cor}
\begin{proof}
After noticing that 
\begin{equation}
\pi^{-1}(gH_1)=\{ ghH: hH\in H_1/H\} ,
\end{equation}
it follows from Lemma \ref{lemsh}  and Theorem \ref{thcov}.
\end{proof} 

\begin{defn}\label{defcovw}
Let $G_1$  and $G_2$ be weighted graphs. 
We'll say that a map $\pi: V(G_1)\rightarrow V(G_2)$ is a left semicovering, iff 
it is surjective  and
for any vertices $u,v\in V(G_2)$, $u_1\in V(G_1)$ such that $u=\pi(u_1)$, we have 
\begin{equation}\label{covw1}
W(G_2)_{uv}=\sum_{v_j\in\pi^{-1}(v)} W(G_1)_{u_1v_j} .
\end{equation} 
We'll say that a map $\pi: V(G_1)\rightarrow V(G_2)$ is a right semicovering, iff
it is surjective  and
for any vertices $u,v\in V(G_2)$, $v_1\in V(G_1)$ such that $v=\pi(v_1)$, we have 
\begin{equation}\label{covw2}
W(G_2)_{uv}=\sum_{u_i\in\pi^{-1}(u)} W(G_1)_{u_iv_1} .
\end{equation}
We'll call a map $\pi: V(G_1)\rightarrow V(G_2)$ a weak covering, iff it is both 
left  and right semicovering. 
\end{defn}

\begin{defn}\label{defcovw1}
Let $G_1$  and $G_2$ be graphs. We'll say that a map $\pi: V(G_1)\rightarrow V(G_2)$ 
is a left semicovering, a right semicovering, or a weak covering if it is a left 
semicovering, a right semicovering, or a weak covering, correspondingly,
of weighted graphs obtained 
from $G_1$  and $G_2$ by choosing $Z=\ZZ$  and $w\equiv 1$.
\end{defn}

Note that left  and right semicoverings as well as a weak covering, are not 
(weighted) graph morphisms. 

\begin{lemma}
For a covering $\pi: G_1\rightarrow G_2$, its vertex map $\pi_V: V(G_1)\rightarrow V(G_2)$ 
is a weak covering.
\end{lemma}
\begin{proof} 
$\pi_V$ is surjective, because $\pi$ is surjective. 
Let $u,v\in V(G_2)$  and $u_1,v_1\in V(G_1)$ be vertices such that 
\begin{equation}
\pi_V(u_1)=u,\quad \pi_V(v_1)=v .
\end{equation}
For each edge $e$ from $u$ to $v$ there is exactly one edge $e_1$ from $u_1$ such that 
$e=\pi_E(e_1)$ and it has the same weight as $e$, because $\pi$ is a local isomorphism. 
For each edge $e_1$ from $u$ to a vertex $v_i\in \pi^{-1}(v)$, the edge $\pi(e_1)$ is 
an edge from $u$ to $v$. Hence, we have the same count of edges from $u$ to $v$  and 
from $u$ to $\pi^{-1}(v)$ and they have the same weights. \eqref{covw1} follows from 
that. We can check \eqref{covw2} analogously, or it follows from \eqref{covw1} by 
changing the directions of all edges.
\end{proof}

\begin{theorem}\label{thcovw}
Let $G_1$  and $G_2$ be (weighted) graphs and let 
$\pi: V(G_1)\rightarrow V(G_2)$ be a map and $\pi(u_1)=u$, $\pi(v_1)=v$ for 
some vertices of $G_1$  and $G_2$. If $\pi$ is a left semicovering, then
\begin{equation}\label{thpiw}
c(G_2,u\rightarrow v)=\sum_{v_j\in\pi^{-1}(v)} c(G_1,u_1\rightarrow v_j) .
\end{equation}
If $\pi$ is a right semicovering, then
\begin{equation}\label{thpiw1}
c(G_2,u\rightarrow v)=\sum_{u_i\in\pi^{-1}(u)}c(G_1,u_i\rightarrow v_1).
\end{equation}
If $\pi$ is a weak covering, then both \eqref{thpiw}  and \eqref{thpiw1} are true.
\end{theorem}
\begin{proof}
Let graphs be weighted  and $\pi$ be a left semicovering.
We'll prove that
\begin{equation}\label{thpiw2}
c_n(G_2,u\rightarrow v)=\sum_{v_j\in\pi^{-1}(v)} c_n(G_1,u_1\rightarrow v_j) 
\end{equation}
is true for all $n$, by induction on $n$. For $n=0$, if $u\neq v$, then $u_1\neq v_j$ for 
all $v_j\in\pi^{-1}(v)$ and both sides of 
\eqref{thpiw2} are $0$;  and if $u=v$,then $u_1\in\pi^{-1}(v)$ and both sides of 
\eqref{thpiw2} are equal to $1$. For $n=1$, \eqref{thpiw2} is equivalent to \eqref{covw1}, 
because 
\begin{equation}\label{thpiw3}
c_1(G, u\rightarrow v)=W(G)_{uv}
\end{equation}
for any $G$, $u$ and $v$. Now let \eqref{thpiw2} be true for some $n\geq 1$ 
for all the correct data. By Definition \ref{def3}, 
\begin{equation}\label{thpiw4}
c_{n+1}(G,u\rightarrow v)=\sum_{x\in V(G)}c_n(G, u\rightarrow x)c_1(G, x\rightarrow v) .
\end{equation} 
for any $G$, $u$, $v$  and $n\geq 0$. For $G=G_2$, substituting \eqref{thpiw2} into the 
right hand side of \eqref{thpiw4}, we get
\begin{gather}\notag
c_{n+1}(G_2,u\rightarrow v)= 
\sum_{x\in V(G_2)}\sum_{x_i\in\pi^{-1}(x)}c_n(G_1, u_1\rightarrow x_i)\sum_{v_j\in\pi^{-1}
(v)} c_1(G_1, x_i\rightarrow v_j)=\\ \label{thpiw5}
\sum_{x_i\in V(G_1)}c_n(G_1, u_1\rightarrow x_i)\sum_{v_j\in\pi^{-1}
(v)} c_1(G_1, x_i\rightarrow v_j)=\\ \notag
\sum_{v_j\in\pi^{-1}(v)}\sum_{x_i\in V(G_1)}c_n(G_1, u_1\rightarrow x_i)
c_1(G_1, x_i\rightarrow v_j)=
\sum_{v_j\in\pi^{-1}(v)}c_{n+1}(G_1, u_1\rightarrow v_j) .
\end{gather}
By induction, \eqref{thpiw2} is true for all $n$, that is equivalent to \eqref{thpiw}.
If $\pi$ is a right semicovering, we can prove \eqref{thpiw1} analogously, or it follows 
from \eqref{thpiw} by changing the directions of all edges. If $\pi$ is a weak covering, 
it is both left  and right covering, so both \eqref{thpiw}  and \eqref{thpiw1} are true.
By Definition \ref{defcovw1}, we can consider weight function $w\equiv 1\in\ZZ$ for
nonweighted graphs.
\end{proof} 

\begin{lemma}\label{lemmap}
If $G_1$ is a positive weighted graph, $G_1$  and $G_2$ are both either directed or 
undirected graphs and $\pi:V(G_1)\rightarrow V(G_2)$ is either a left 
semicovering, or a right semicovering, or a weak covering, then there is such map 
$\pi_E: E(G_1)\rightarrow E(G_2)$ that the pair $(\pi, \pi_E)$ is a graph morphism from 
$G_1$ to $G_2$.
\end{lemma}
\begin{proof}
Let $u=\pi(u_1)$  and $v=\pi(v_1)$.
For each edge $e_1\in E(G_1)$ from $u_1$ to $v_1$ we have 
\begin{equation}
W_{uv}=W_{u_1v_1}+\dots=w(e_1)+\dots>0 .
\end{equation}
It means that there is an edge $e$ from $u$ to $v$ and we can put $\pi_E(e_1)=e$. 
\end{proof}

Note that we constructed only a graph morphism, not a weighted graph morphism. It is possible 
that the weighted graph morphism continuing a semicovering, or a weak covering, doesn't 
exist.

\begin{lemma}\label{lemmodw1}
The map $\pi: V(|G_1|)\rightarrow V(|G_2|)$ is a left semicovering, a right semicovering, or 
a weak covering, iff it is a left semicovering, a right semicovering, or a weak covering, 
correspondingly, considered as a map $\pi: V(G_1)\rightarrow V(G_2)$.
\end{lemma}
\begin{proof}
Weight matrix elements are the same for $G$  and $|G|$ and Definition \ref{defcovw} depends 
only on the matrix elements.  
\end{proof}  

\begin{cor}\label{corpos}
If $G_1$ is a connected positive weighted graph and $\pi: V(G_1)\rightarrow G_2$ is either 
a left semicovering, or a right semicovering, or a weak covering, then $G_2$ is a connected 
weighted graph.
\end{cor}
\begin{proof}
By Lemma \ref{lemmodw1}, we can replace $G_1$  and $G_2$ by undirected graphs $|G_1|$  and 
$|G_2|$. By Lemma \ref{lemmap}, the map $\pi$ can be continued to a graph morphism from 
$|G_1|$ to $|G_2|$. Now, $\pi$ is surjective and the image of a connected graph under 
a graph morphism is connected. 
\end{proof} 

\begin{cor}
If $G_1$ is a positive weighted graph, $G_2$ is a reduced weighted graph and 
$\pi: V(G_1)\rightarrow V(G_2)$ is either a left semicovering, or a right semicovering, or 
a weak covering, then $G_2$ is a positive weighted graph.
\end{cor}
\begin{proof}
If $G_2$ is reduced, then the weights
of its edges coincide with the corresponding matrix elements of its weight matrix and by 
Definition \ref{defcovw} they are sums of positive elements, hence positive. 
\end{proof}

\begin{cor}
If $G_1$ is a positive weighted graph, $G_2$ is a reduced weighted graph and 
$\pi: V(G_1)\rightarrow V(G_2)$ is a weak covering, then   
for each connected component $G$ of $G_1$, 
the complete subgraph with the set of vertices $\pi(V(G))\subseteq V(G_2)$ is a connected
component of $G_2$.
\end{cor}
\begin{proof}
Let $u_1\in G$  and $u=\pi(u_1)$. For each edge from $u$, we can find an edge from $u_1$ 
projecting onto it, because $\pi$ is a left semicovering and the corresponding matrix 
element of the weight matrix of $G_2$ is nonzero. End of the edge starting in $u_1$, lies 
in $G$, because $G$ is a connected component. Hence, the projection of this end lies in 
$\pi(G)$. Analogously for edges ending in $u$ we can prove that they start in $\pi(G)$, 
using that $\pi$ is a right semicovering. 
\end{proof}

\begin{theorem}\label{thcovw1}
Let $G_1$ be a connected positive weighted graph  and $\pi: V(G_1)\rightarrow V(G_2)$ 
be a weak covering. Then $G_2$ is a connected weighted graph as well and either the sets 
$\pi^{-1}(v)$ are infinite for each vertex $v\in \pi(G)$, or all of them are finite 
 and contain an equal number of elements.
\end{theorem}
\begin{proof}
The connectivity of $G_2$ is already stated in Corollary \ref{corpos}.
The other proof is the same as for Theorem \ref{thcov1}, except that we consider 
weighted walks now and we need a new explanation why we can change the order of summation 
in \eqref{eqcov2}. The point is that the sum have only a finite number of coefficients 
at $t^n$ for each $n$. Also, instead of \eqref{eqcov4} we need 
\begin{equation}\label{eqcovw4}
c(G_2, u\rightarrow v)> 0 
\end{equation}
that is true because $c(G_2, u\rightarrow v)$ is a sum of positive summands.
\end{proof}

\begin{cor}\label{corcovw}
Let $G_1$ be a positive weighted graph, 
 and $\pi: V(G_1)\rightarrow V(G_2)$ be a weak covering. 
Then for each connected component $G$ of $G_1$ either the sets 
$\pi^{-1}(v)\cap G$ are infinite for each vertex $v\in \pi(G)$, or all of them are finite 
 and contain an equal number of elements.
\end{cor}
\begin{proof}
The restriction $\pi|_G:V(G)\rightarrow \pi(G)$ is a weak covering from $G$ to the complete 
subgraph of $G_2$ with the set of vertices $\pi(G)$ and we can use Theorem \ref{thcovw1}.
\end{proof}   
 
\begin{defn}
Denote $K_1(m)=(0, 0, m\cdot (0))$ 
the graph with one vertex  and one loop of weight $m$. 
\end{defn}

\begin{theorem}\label{thKm}
$c(K_1(m), 0)=e^{mt}$ .
\end{theorem}
\begin{proof}
The weight matrix of $K_1(m)$ is $(m)$ and we can use Lemma \ref{lem2}.
\end{proof}

\begin{cor}\label{corKm}
If a weighted graph $G_1$ is obtained from a weighted graph $G$ by adding loops of weight $m$
to all the vertices, then 
\begin{equation}\label{eqcorKm}
c(G_1,u\rightarrow v)=e^{mt}c(G,u\rightarrow v)
\end{equation}
for any their vertices $u, v$.
\end{cor}
\begin{proof}
$G_1$ is isomorphic to $G\times K_1^m$ and we can use Theorems \ref{thm1}  and \ref{thKm}.
\end{proof}

\begin{defn}\label{defreg}
We'll call a weighted graph $G$ left semiregular of weight $m$ iff for each vertex $u$ 
\begin{equation}
\sum_{v\in V(G)}W(G)_{uv}=m .
\end{equation}
right semiregular of weight $m$ iff for each vertex $v$ 
\begin{equation}
\sum_{u\in V(G)}W(G)_{uv}=m , 
\end{equation}
 and weak regular of weight $m$ if it is both left semiregular of weight $m$  and 
right semiregular of weight $m$.
\end{defn}

\begin{theorem}\label{threg}
A weighted graph $G$ is left semiregular of weight $m$, or right semiregular of weight $m$, or 
weak regular of weight $m$ iff the constant map $\pi: V(G)\rightarrow V(K_1(m))$ is a left 
semicovering, or a right semicovering, or a weak covering, correspondingly.
\end{theorem}
\begin{proof}
It is enough to compare Definitions \ref{defcovw}  and \ref{defreg}. 
\end{proof} 

\begin{cor}\label{correg}
For a left semiregular graph $G$ of degree $m$ and a vertex $u\in V(G)$,
\begin{gather}\label{reg1}
\sum_{v\in V(G)}c(G, u\rightarrow v)=e^{mt},\\
\sum_{v\in V(G)}c_n(G,u\rightarrow v)=m^n.
\label{reg2}
\end{gather}
For a right semiregular graph $G$ of degree $m$ and a vertex $v\in V(G)$,
\begin{gather}\label{reg3}
\sum_{u\in V(G)}c(G, u\rightarrow v)=e^{mt},\\
\sum_{u\in V(G)}c_n(G,u\rightarrow v)=m^n.
\label{reg4}
\end{gather}
For a weak regular graph $G$ of degree $m$ all the formulas $(\ref{reg1}-\ref{reg4})$ 
are true for any vertices $u,v\in V(G)$.
\end{cor}
\begin{proof}
It is an immediate corollary of Theorems \ref{thcovw}, \ref{thKm} and  \ref{threg}.
\end{proof}

\section{Symmetric powers  and bipowers}

\begin{defn}\label{nd1}
For a (weighted) graph $G$ with linearly ordered set of vertices 
$V(G)$, denote $S^n(G)$ a directed (weighted) graph with the set of vertices 
\begin{equation}\label{d3eq1} 
V(S^nG)=S^nV(G)=\{(v_1,\dots,v_n)\in V(G)^n\quad |\quad v_1\leq\dots\leq v_n\} 
\end{equation}
 and edges defined as follows. 
For $i=1,\dots,n$ denote $\Pr_i:V(S^nG)\rightarrow V(S^{n-1}G)$ the projection 
obtained by the erasing of the vertices $v_i$ standing on the $i$-th place of the 
sequence \eqref{d3eq1}. The directed edges from $A$ to $B$ correspond to pairs 
$(i\in \{1,\dots,n\},\thinspace e\in E(G))$ such that    
there is $j$ such that 
$Pr_i(A)=Pr_j(B)$  and $e$ is an edge from 
$v_i(A)$ to $v_j(B)$ where $v_i$ denotes the vertex of $G$ standing on the 
$i$-th place of the sequence \eqref{d3eq1}; in which case, for weighted graphs,
the corresponding edge from $A$ to $B$ has the same weight as $e$. 
\end{defn} 

\begin{lemma}\label{nl1}
For a (weighted) graph $G$ with linearly ordered set of vertices, mapping 
\begin{equation}\label{nl1eq1}
V(G^n)\rightarrow V(S^nG),\quad (v_1,\dots,v_n)\mapsto (v_1,\dots,v_n)^\sigma
\end{equation}
where the permutation $\sigma$ rearranges $(v_1,\dots,v_n)$ in increasing order, 
is a left semicovering. 
\end{lemma}
\begin{proof} 
It follows directly from Definitions \ref{nd1}  and \ref{defcovw}.
\end{proof}

\begin{theorem}\label{nt1}
For a (weighted) graph $G$ with linearly ordered set of vertices, 
\begin{equation}\label{nt1eq1}
c(S^nG, (u_1,\dots,u_n)\rightarrow (v_1,\dots,v_n) )=\per(c(G,u_i\rightarrow v_j))_{1\leq i,j\leq n} , 
\end{equation}
i.\ e.\ the permanent of the matrix with pointed out elements.
\end{theorem}
\begin{proof}
It follows from Lemma \ref{nl1}  and Theorems \ref{thcovw}  and \ref{thm1}, that  
\begin{gather}\label{nt1eq2}
c(S^nG, (u_1,\dots,u_n)\rightarrow (v_1,\dots,v_n) )
=\sum_{\sigma\in S_n} c(G^n,(u_1,\dots,u_n)\rightarrow (v_1,\dots,v_n)^\sigma)\\ 
=\sum_{\sigma\in S_n} c(G, u_1\rightarrow v_{\sigma^{-1}(1)}) \dots 
c(G, u_n\rightarrow v_{\sigma^{-1}(n)})=\per(c(G,u_i\rightarrow v_j))_{1\leq i,j\leq n} .
\label{nt1eq3}
\end{gather}
\end{proof}

\begin{defn}\label{nd2}
For a (weighted) directed graph $G$ with linearly ordered set of vertices 
$V(G)$, denote $S^n_2(G)$ a directed (weighted) graph with the set of vertices 
\begin{equation}\label{nd2eq1} 
V(S^n_2G)=S^nV(G)=\{(v_1,\dots,v_n)\in V(G)^n\quad |\quad v_1\leq\dots\leq v_n\} 
\end{equation}
 and edges defined as follows. The directed edges from $A$ to $B$ correspond to 
the pairs $(\sigma\in S_n, (e_1,\dots,e_n)\in (E(G))^n)$ such that 
for $i$ from $1$ to $n$, $e_i$ is an edge from $v_i(A)$ to $v_{\sigma(j)}(B)$ 
where $v_i$ denotes the vertex of $G$ standing on the 
$i$-th place of the sequence \eqref{nd2eq1}, 
in which case, for weighted graphs,
the corresponding edge from $A$ to $B$ has the weight $e_1\dots e_n$.
\end{defn}

\begin{lemma}\label{nl2}
For a (weighted) directed graph $G$ with linearly ordered set of vertices, mapping 
\begin{equation}\label{nl2eq1}
V(\underbrace{G\times_2\dots\times_2 G}_{n\thickspace \text{times}} )
\rightarrow V(S^n_2G),\quad (v_1,\dots,v_n)\mapsto (v_1,\dots,v_n)^\sigma
\end{equation}
where the permutation $\sigma$ rearranges $(v_1,\dots,v_n)$ in increasing order, 
is a left semicovering. 
\end{lemma}
\begin{proof} 
It follows directly from Definitions \ref{nd2}  and \ref{defcovw}.
\end{proof}

\begin{theorem}\label{nt2}
For a (weighted) directed graph $G$ with linearly ordered set of vertices, 
\begin{equation}\label{nt2eq1}
c(S^n_2G, (u_1,\dots,u_n)\rightarrow (v_1,\dots,v_n) )=\per^\circ (c(G,u_i\rightarrow v_j))_{1\leq i,j\leq n} , 
\end{equation}
i.\ e.\ the permanent of the matrix with pointed out elements multiplying 
coefficient-wise according to Definition \ref{dcirc}.
\end{theorem}
\begin{proof}
It follows from Lemma \ref{nl2}  and Theorems \ref{thcovw}  and \ref{thbi1}, 
the same as in the proof of Theorem \ref{nt2}.
\end{proof}

\begin{cor}\label{nc1}
For a (weighted) directed graph $G$ with linearly ordered set of vertices, 
\begin{equation}\label{nc1eq1}
c_k(S^n_2G, (u_1,\dots,u_n)\rightarrow (v_1,\dots,v_n) )=\per (c_k(G,u_i\rightarrow v_j))_{1\leq i,j\leq n} , 
\end{equation}
i.\ e.\ the permanent of the matrix with pointed out elements. 
\end{cor}
\begin{proof}
It follows directly from Theorem \ref{nt2}  and Definition \ref{dcirc}.
\end{proof}

\begin{defn}\label{nd3}
For a bipartite graph $G$, we'll say that a vertex $(v_1,\dots,v_n) \in V(S^n_2G)$ 
is even if the sum of parities $(0,1)$ of $v_1, \dots, v_n$ is even. If the vertex 
is not even, we'll say that it is odd.
\end{defn}

\begin{lemma}\label{nl3}
Let $G$ be a bipartite directed graph. 
For even $n$, each edge of $S^n_2G$ connects vertices of the same parity, 
i.\ e.\ either 
both even, or both odd. For odd $n$, each edge of $S^n_2G$ connects vertices of  
alternate parities, 
i.\ e.\ one of them must be even  and another odd.  
\end{lemma}
\begin{proof}
For an edge from $(u_1,\dots,u_n)$ to $(v_1,\dots,v_n)$, the sum of parities of 
all the $u_i$  and $v_i$ equals $n$, therefore it is even for even $n$  and odd for 
odd $n$.
\end{proof}

\begin{defn}\label{nd4}
For a bipartite directed graph $G$  and even $n$ denote $S^n_0G$  and $S^n_1G$ 
the complete subgraphs of 
$S^n_2G$, spanned by all the even vertices, or the odd, respectively.
\end{defn}

\begin{cor}\label{nc2}
For even $n$, we can use the formulas \eqref{nt2eq1}  and \eqref{nc1eq1} for 
$S^n_0G$  and $S^n_1G$ instead of $S^n_2G$, as well.
\end{cor}
\begin{proof}
It follows from Theorem \ref{nt2}, Corollary \ref{nc1}  and Lemma \ref{nl3}.
\end{proof}

\section{Schur operations  and bioperations}

We'll fix here for each Young diagram of a partition $\lambda$, a standard tableau $T$ of shape 
$\lambda$ filling first the first row, then the second  and so on. 

\begin{defn}\label{nd6}
For a (weighted) graph $G$ with linearly ordered set of vertices 
$V(G)$  and for a partition $\lambda$ denote $S^\lambda(G)$ 
a directed (weighted) graph vertices 
of which are tableaux of shape 
$\lambda$ with elements from $V(G)$ increasing in each column  and not decreasing 
in each row,  
 and edges defined as follows. 
For a triple consisting of a vertex $u\in V(S^\lambda G)$, box $i$ of the 
Young diagram of shape $\lambda$  and an edge $e$ from the vertex $u_i\in V(G)$ 
located on the $i$-th place of the Young diagram, to a vertex $v\in V(G)$, 
denote $(u,i,e)$ the tableau obtained by substituting $v$ instead of $u_i$ in 
the $i$-th box of $u$. It is not a tableau with increasing columns  and non-decreasing 
rows, in general. Let  
\begin{equation}\label{nd6eq1}
Y_T(u,i,e)=\sum_{x\in V(S^\lambda G)} \gamma_x Y_T(x)
\end{equation} 
where $Y_T$ is the Young symmetrizer corresponding to the standard tableau 
$T$ fixed in the beginning of this section. $\gamma_x\in \ZZ$ are some integer 
coefficients. We suppose that for each such triple $(u,i,v)$  and for each 
non-zero coefficient $\gamma_x$ in \eqref{nd6eq1} we have an edge from $u$ to $x$ 
of weight $\gamma_x w(e)$ where $w(e)$ is the weight of $e$.
\end{defn}

\begin{theorem}\label{nt3}
For a (weighted) graph $G$ with linearly ordered set of vertices  and for a partition 
$\lambda$, 
\begin{equation}\label{nt3eq1}
c(S^\lambda G, (u_1,\dots,u_n)\rightarrow (v_1,\dots,v_n) )=\imm_\lambda(c(G,u_i\rightarrow v_j))_{1\leq i,j\leq n} , 
\end{equation}
i.\ e.\ the immanent of the matrix with pointed out elements, where 
\begin{equation}\label{nt3eq2}
\imm_\lambda(a_{ij})_{1\leq i,j\leq n}=\sum_{\sigma\in S_n}\chi_\lambda a_{1\sigma(1)}\dots a_{n\sigma(n)}
\end{equation}
 and $\chi_\lambda$ denotes the character of symmetric group $S_n$ corresponding to 
the partition $\lambda$. 
\end{theorem}
\begin{proof}
It follows from the commutativity of the diagram 
\begin{equation}\label{nt3eq3}
\begin{CD}
\Mat_{V(G)}(r) @>\exp>>\Mat_{V(G)}(r\{\{t\}\})\\
@V{S^\lambda}VV @VV{S^\lambda_\ast}V\\
\Mat_{V(S^\lambda G)}(r) @>\exp>>\Mat_{V(S^\lambda G)}(r\{\{t\}\}) 
\end{CD}
\end{equation}
for our ring of weights $r$, 
where $\exp: W(G)\mapsto C(G)=e^{W(G)t}$, $S^\lambda (W(G))=W(S^\lambda G)$ 
 and $S^\lambda_\ast$ is defined according to the right hand side of \eqref{nt3eq1}. 
The commutativity of this diagram for $r$ follows from the commutativity of 
this diagram for $r=\ZZ$ which follows from the commutativity of this 
diagram for $r=\RR$ in which case it follows from the well-known correspondence 
between the representations of Lie groups  and algebras. Also, $S^\lambda$, 
$S^\lambda_\ast$  and $\exp$ can be defined as univeral functors for appropriate 
categories. Combining the universal properties, we obtain a unique universal 
diagonal functor for \eqref{nt3eq3}.
\end{proof}

Note that Theorem \ref{nt1} is a particular case of Theorem \ref{nt3} for 
$\lambda=n$.

\begin{defn}\label{nd7}
For a (weighted) graph $G$ with linearly ordered set of vertices 
$V(G)$  and for a partition $\lambda$ denote $S^\lambda_2(G)$ 
a directed (weighted) graph vertices 
of which are tableaux of shape 
$\lambda$ with elements from $V(G)$ increasing in each column  and not decreasing 
in each row,  
 and edges defined as follows. 
For a pair consisting of a vertex $u\in V(S^\lambda G)$  and tableau $e$  
of shape $\lambda$ filled by the edges of $G$ such that for each box $i$ the 
edge $e_i$ is an edge from the vertex $u_i\in V(G)$ 
located in the $i$-th box of $u$ to a vertex $v_i\in V(G)$, 
denote $(e)$ the tableau obtained by substituting $v_i$ instead of $u_i$ in 
the $i$-th box of $u$. It is not a tableau with increasing columns  and non-decreasing 
rows, in general. Let  
\begin{equation}\label{nd7eq1}
Y_T(e)=\sum_{x\in V(S^\lambda G)} \gamma_x Y_T(x)
\end{equation} 
where $Y_T$ is the Young symmetrizer corresponding to the standard tableau 
$T$ fixed in the beginning of this section. $\gamma_x\in \ZZ$ are some integer 
coefficients. We suppose that for each an edge tableau $(e)$  and for each 
non-zero coefficient $\gamma_x$ in \eqref{nd6eq1} we have an edge from $u$ to $x$ 
of weight $\gamma_x w(e_1)\dots w(e_n)$ where $w(e_i)$ is the weight of $e_i$.
\end{defn}

\begin{theorem}\label{nt4}
For a (weighted) graph $G$ with linearly ordered set of vertices  and for a partition 
$\lambda$, 
\begin{equation}\label{nt4eq1}
c(S^\lambda_2 G, (u_1,\dots,u_n)\rightarrow (v_1,\dots,v_n) )=\imm_\lambda^\circ(c(G,u_i\rightarrow v_j))_{1\leq i,j\leq n} , 
\end{equation}
i.\ e.\ the immanent of the matrix with pointed out elements 
multiplying 
coefficient-wise according to Definition \ref{dcirc}, i.\ e.\  
\begin{equation}\label{nt4eq2}
\imm_\lambda^\circ(a_{ij})_{1\leq i,j\leq n}=
\sum_{\sigma\in S_n}\chi_\lambda a_{1\sigma(1)}\circ\dots\circ a_{n\sigma(n)}
\end{equation}
 and $\chi_\lambda$ denotes the character of symmetric group $S_n$ corresponding to 
the partition $\lambda$. 
\end{theorem}
\begin{proof}
Analogously the proof of Theorem \ref{nt3}, it follows 
from the commutativity of the diagram 
\begin{equation}\label{nt4eq3}
\begin{CD}
\Mat_{V(G)}(r) @>\exp>>\Mat_{V(G)}(r\{\{t\}\})\\
@V{S^\lambda_2}VV @VV{S^\lambda_\circ}V\\
\Mat_{V(S^\lambda_2 G)}(r) @>\exp>>\Mat_{V(S^\lambda_2 G)}(r\{\{t\}\}) 
\end{CD}
\end{equation}
for our ring of weights $r$, 
where $\exp: W(G)\mapsto C(G)=e^{W(G)t}$, $S^\lambda(W(G))=W(S^\lambda G)$ 
 and $S^\lambda_\circ$ is defined according to the right hand side of \eqref{nt3eq1}. 
Again, $S^\lambda_2$, 
$S^\lambda_\circ$  and $\exp$ can be defined as univeral functors for appropriate 
categories. Combining the universal properties, we obtain a unique universal 
diagonal functor for \eqref{nt4eq3}.
\end{proof}

\begin{cor}\label{nc3}
For a (weighted) graph $G$ with linearly ordered set of vertices  and for a partition 
$\lambda$, 
\begin{equation}\label{nc3eq1}
c_k(S^\lambda_2 G, (u_1,\dots,u_n)\rightarrow (v_1,\dots,v_n) )=\imm_\lambda(c(G,u_i\rightarrow v_j))_{1\leq i,j\leq n} , 
\end{equation}
i.\ e.\ the immanent of the matrix with pointed out elements. 
\end{cor}
\begin{proof}
It follows directly from Theorem \ref{nt4}  and Definition \ref{dcirc}.
\end{proof}

Note that Theorem \ref{nt1}  and Corollary \ref{nt1} are the particular cases of 
Theorem \ref{nt4}  and Corollary \ref{nc3} for $\lambda=n$.

\section{Exterior powers  and bipowers}   

\begin{defn}\label{nd5}
For a (weighted) graph $G$ with linearly ordered set of vertices 
$V(G)$, denote $\Lambda^n(G)$ a directed (weighted) graph with the set of vertices 
\begin{equation}\label{d5eq1} 
V(\Lambda^nG)=\Lambda^nV(G)=\{(v_1,\dots,v_n)\in V(G)^n\quad |\quad v_1 < \dots < v_n\} 
\end{equation}
 and edges defined as follows. 
For $i=1,\dots,n$ denote $\Pr_i:V(\Lambda^nG)\rightarrow V(\Lambda^{n-1}G)$ the projection 
obtained by the erasing of the vertices $v_i$ standing on the $i$-th place of the 
sequence \eqref{d5eq1}. The directed edges from $A$ to $B$ correspond to pairs 
$(i\in \{1,\dots,n\},\thinspace e\in E(G))$ such that    
there is $j$ such that 
$Pr_i(A)=Pr_j(B)$  and $e$ is an edge from 
$v_i(A)$ to $v_j(B)$ where $v_i$ denotes the vertex of $G$ standing on the 
$i$-th place of the sequence \eqref{d5eq1}; in which case, for weighted graphs,
the corresponding edge from $A$ to $B$ has the weight $(-1)^{j-i}w(e)$ where 
$w(e)$ is the weight of $e$. 
\end{defn} 

\begin{theorem}\label{nt5}
For a (weighted) graph $G$ with linearly ordered set of vertices, 
\begin{equation}\label{nt5eq1}
c(\Lambda^nG, (u_1,\dots,u_n)\rightarrow (v_1,\dots,v_n) )=
\det(c(G,u_i\rightarrow v_j))_{1\leq i,j\leq n} . 
\end{equation}
\end{theorem}
\begin{proof}
It is a particular case of Theorem \ref{nt3} for $\lambda=1^n$.
\end{proof}

\begin{defn}\label{nd8}
For a (weighted) directed graph $G$ with linearly ordered set of vertices 
$V(G)$, denote $\Lambda^n_2(G)$ a directed (weighted) graph with the set of vertices 
\begin{equation}\label{nd8eq1} 
V(\Lambda^n_2G)=\Lambda^nV(G)=\{(v_1,\dots,v_n)\in V(G)^n\quad |\quad v_1\leq\dots\leq v_n\} 
\end{equation}
 and edges defined as follows. The directed edges from $A$ to $B$ correspond to 
the pairs $(\sigma\in S_n, (e_1,\dots,e_n)\in (E(G))^n)$ such that 
for $i$ from $1$ to $n$, $e_i$ is an edge from $v_i(A)$ to $v_{\sigma(j)}(B)$ 
where $v_i$ denotes the vertex of $G$ standing on the 
$i$-th place of the sequence \eqref{nd8eq1}, 
in which case, for weighted graphs,
the corresponding edge from $A$ to $B$ has the weight $\sgn(\sigma)e_1\dots e_n$.
\end{defn}

\begin{theorem}\label{nt6}
For a (weighted) directed graph $G$ with linearly ordered set of vertices, 
\begin{equation}\label{nt6eq1}
c(\Lambda^n_2G, (u_1,\dots,u_n)\rightarrow (v_1,\dots,v_n) )=
\det^\circ (c(G,u_i\rightarrow v_j))_{1\leq i,j\leq n} , 
\end{equation}
i.\ e.\ the determinant of the matrix with pointed out elements multiplying 
coefficient-wise according to Definition \ref{dcirc}.
\end{theorem}
\begin{proof}
It is a particular case of Theorem \ref{nt4} for $\lambda=1^n$.
\end{proof}

\begin{cor}\label{nc4}
For a (weighted) directed graph $G$ with linearly ordered set of vertices, 
\begin{equation}\label{nc4eq1}
c_k(\Lambda^n_2G, (u_1,\dots,u_n)\rightarrow (v_1,\dots,v_n) )=\det (c_k(G,u_i\rightarrow v_j))_{1\leq i,j\leq n} , 
\end{equation}
i.\ e.\ the determinant of the matrix with pointed out elements. 
\end{cor}
\begin{proof}
It follows directly from Theorem \ref{nt6}  and Definition \ref{dcirc}.
\end{proof}

\begin{defn}\label{nd9}
For a bipartite graph $G$, we'll say that a vertex $(v_1,\dots,v_n) \in V(\Lambda^n_2G)$ 
is even if the sum of parities $(0,1)$ of $v_1, \dots, v_n$ is even. If the vertex 
is not even, we'll say that it is odd.
\end{defn}

\begin{lemma}\label{nl4}
Let $G$ be a bipartite directed graph. 
For even $n$, each edge of $\Lambda^n_2G$ connects vertices of the same parity, 
i.\ e.\ either 
both even, or both odd. For odd $n$, each edge of $\Lambda^n_2G$ connects vertices of  
alternate parities, 
i.\ e.\ one of them must be even  and another odd.  
\end{lemma}
\begin{proof}
For an edge from $(u_1,\dots,u_n)$ to $(v_1,\dots,v_n)$, the sum of parities of 
all the $u_i$  and $v_i$ equals $n$, therefore it is even for even $n$  and odd for 
odd $n$.
\end{proof}

\begin{defn}\label{nd10}
For a bipartite directed graph $G$  and even $n$ denote $\Lambda^n_0G$  and 
$\Lambda^n_1G$ the complete subgraphs of 
$\Lambda^n_2G$, spanned by all the even vertices, or the odd, respectively.
\end{defn}

\begin{cor}\label{nc5}
For even $n$, we can use the formulas \eqref{nt6eq1}  and \eqref{nc4eq1} for 
$\Lambda^n_0G$  and $\Lambda^n_1G$ instead of $\Lambda^n_2G$, as well.
\end{cor}
\begin{proof}
It follows from Theorem \ref{nt6}, Corollary \ref{nc4}  and Lemma \ref{nl4}.
\end{proof}

\section{Linear graphs}\label{sec2}

\begin{defn}
Denote $R=(\ZZ, 0, 1, -1)$ the  graph with the set of 
vertices $\ZZ$ and edges connecting $n$  and $n+1$ for every n:

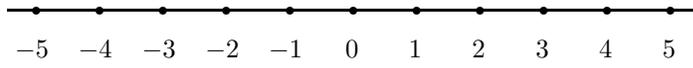
\begin{picture}(410, 25) 
\put(83,18){\circle*{3}}
\put(107,18){\circle*{3}}
\put(131,18){\circle*{3}}
\put(155,18){\circle*{3}}
\put(179,18){\circle*{3}}
\put(203,18){\circle*{3}}
\put(227,18){\circle*{3}}
\put(251,18){\circle*{3}}
\put(275,18){\circle*{3}}
\put(299,18){\circle*{3}}
\put(323,18){\circle*{3}}
\put(75,0){$-5$}
\put(99,0){$-4$}
\put(123,0){$-3$}
\put(147,0){$-2$}
\put(171,0){$-1$}
\put(200,0){$0$}
\put(224,0){$1$}
\put(248,0){$2$}
\put(272,0){$3$}
\put(296,0){$4$}
\put(320,0){$5$}
\put(84,18){\line(1,0){22}}
\put(108,18){\line(1,0){22}}
\put(132,18){\line(1,0){22}}
\put(156,18){\line(1,0){22}}
\put(180,18){\line(1,0){22}}
\put(204,18){\line(1,0){22}}
\put(228,18){\line(1,0){22}}
\put(252,18){\line(1,0){22}}
\put(276,18){\line(1,0){22}}
\put(300,18){\line(1,0){22}}
\put(72,18){\line(1,0){10}}
\put(324,18){\line(1,0){10}}
\end{picture}
\end{defn}

\begin{theorem}\label{thm3}
\begin{equation}\label{19}
c_n(R,m)=\begin{cases}\dbinom{n}{\frac{n-m}{2}} &\text{if $n-m$ is even,}\\
0 &\text{otherwise.} \end{cases}
\end{equation}
\end{theorem}
\begin{proof}
Denoting $q$ the element of the group algebra $\ZZ(\ZZ)$, corresponding to $1\in\ZZ$ (inside 
parentheses), we get
\begin{equation}\label{19a}
\ZZ(\ZZ)=\ZZ[q,q^{-1}] ,
\end{equation}
where $q^m$ corresponds to $m$.
Then by \eqref{eqcorexp1},
\begin{equation}\label{19b}
(q+q^{-1})^n=\sum_{m\in\ZZ}c_n(R,m)q^m .
\end{equation}
 and we know that 
\begin{equation}\label{20}
(q+q^{-1})^n=\sum_{i=0}^n\binom{n}{i}q^{n-2i},
\end{equation}
so
\begin{equation}\label{21}
c_n(R,m)=\binom{n}{i} ,
\end{equation}
where $m=n-2i$, i.\ e.\ $i=(n-m)/2$ .
\end{proof}

\begin{cor}\label{cor4}For $m\geq 0$,
\begin{equation}\label{24}
c(R,\pm m)=\frac{t^{m}}{m!}\thickspace {_0 F_1}(;m+1;t^2)=I_m(2t) ,
\end{equation}
where $_0F_1$ is the hypergeometric function,
 and $I_m$ is the modified Bessel function of the first kind.
\end{cor}
\begin{proof}It follows immediately from Theorem \ref{thm3}  and the definitions of 
hypergeometric functions  and Bessel functions.
\end{proof}

\begin{cor}\label{cor4a}
\begin{equation}\label{eq4a}
e^{(q+q^{-1})t}=\sum_{m=-\infty}^{\infty}q^mI_{|m|}(2t) .
\end{equation}
\end{cor}
\begin{proof} 
It follows from \eqref{eqcorexp2}  and Corollary \ref{cor4}. 
\end{proof}

\begin{defn}Denote $R(a,b)=(\ZZ, 0, a\cdot 1, b\cdot (-1)) $ the graph $R$, 
considered as a directed graph having for every 
integer $n$ the arcs $(n,n+1)$ of weight $a$  and $(n+1,n)$ of weight $b$.
\end{defn} 

\begin{theorem}\label{thm4}
\begin{equation}\label{25}
c_n(R(a,b),m)=\begin{cases}\dbinom{n}{\frac{n-m}{2}}a^{\frac{n+m}{2}}b^{\frac{n-m}{2}}
&\text{if $n-m$ is even,}\\0 &\text{otherwise.} \end{cases}
\end{equation}
\end{theorem}
\begin{proof}
We can 
repeat the proof of Theorem \ref{thm3} with corresponding modifications. Also, we get another 
proof, noticing that each walk of length $n$ from $0$ to $m$ has the same weight $a^{(n+m)/2}
b^{(n-m)/2}$ and using Theorem \ref{thm3}. 
\end{proof}

\begin{cor}\label{cor5}For $m\geq 0$  and $ab\neq 0$,
\begin{gather}
c(R(a,b),m)=\frac{a^mt^m}{m!}\thickspace {_0F_1}(;m+1;abt^2)=a^{m/2}b^{-m/2}
I_m(2\sqrt{ab}\thinspace t) ,\\
c(R(a,b),-m)=\frac{b^mt^m}{m!}\thickspace {_0F_1}(;m+1;abt^2)=a^{-m/2}b^{m/2}
I_m(2\sqrt{ab}\thinspace t) .
\end{gather}
For $ab=0$
\begin{equation}\label{26}
c(R(a,b),m)=\begin{cases}a^mt^m/m! &\text{if $b=0$, $m\geq 0$,}\\
b^{-m}t^{-m}/(-m)! &\text{if $a=0$, $m\leq 0$,}\\
0 &\text{otherwise}.
\end{cases}
\end{equation}
\end{cor}
\begin{proof}
It follows from Theorem \ref{thm4}  and definitions.
\end{proof}

\begin{cor}\label{cor5a}For $ab\neq 0$,
\begin{equation}\label{eq5a}
e^{(aq+bq^{-1})t}=\sum_{m=-\infty}^{\infty}(\sqrt{a/b}\thinspace q)^mI_{|m|}
(2\sqrt{ab}\thinspace t) .
\end{equation}
\end{cor}
\begin{proof} 
The same proof as for Corollary \ref{cor4a}. Or we can change in (\ref{eq4a}) 
$q$ to $\sqrt{a/b}\thinspace q$  and $t$ to $\sqrt{ab}\thinspace t$.
\end{proof}

\begin{theorem}\label{thmth}
\begin{equation}\label{eqthth1}
c_n(R, m_1, \dots, m_k)=c_n(R,m)=
\begin{cases}\dbinom{n}{\frac{n-m}{2}} &\text{if $n-m$ is even,}\\
0 &\text{otherwise,} \end{cases}
\end{equation} 
where $m=|m_1|+\dots+|m_k|$.
\end{theorem}
\begin{proof}
Consider the automorphism 
\begin{equation}
- :R\rightarrow R,\quad m\mapsto (-m) . 
\end{equation} 
Nonnegative integers form a fundamental domain with respect to it and our result 
follows from Theorems \ref{thm3}  and \ref{thfun1}.
\end{proof}  

\begin{cor}\label{corth}
\begin{equation}
c(R,m_1,\dots,m_k)=I_m(2t) ,
\end{equation}
where $m=|m_1|+\dots+|m_k|$.
\end{cor}
\begin{proof}
It follows from Theorem \ref{thmth}  and Corollary \ref{cor4}.
\end{proof}

\begin{theorem}\label{thmthw}
\begin{equation}\label{eqthth1w}
c_n(R(a,b), m_1, \dots, m_k)=
\begin{cases}\dbinom{n}{\frac{n-m}{2}}a^{(n+M)/2}b^{(n-M)/2} &\text{if $n-m$ is even,}\\
0 &\text{otherwise,} \end{cases}
\end{equation} 
where $m=|m_1|+\dots+|m_k|$  and $M=m_1+\dots+m_k$.
\end{theorem}
\begin{proof}
In the proof of Theorem \ref{thm4} we noticed that each walk from $0$ to $M$ 
has the same weight $a^{(n+M)/2}b^{(n-M)/2}$. Multiplying it on the number of walks stated 
in Theorem \ref{thmth}, we obtain (\ref{eqthth1w}).
\end{proof} 

\begin{cor}\label{corthw}For $ab\neq 0$,
\begin{equation}\label{eqcorthw1}
c(R(a,b),m_1,\dots,m_k)=\frac{a^{(m+M)/2}b^{(m-M)/2}t^m}{m!}
\thickspace {_0F_1}(;m+1;abt^2)=a^{M/2}b^{-M/2}
I_m(2\sqrt{ab}\thinspace t) ,
\end{equation}
where $m=|m_1|+\dots+|m_k|$  and $M=m_1+\dots+m_k$.
For $ab=0$
\begin{equation}\label{eqcorthw2}
c(R(a,b),m_1,\dots,m_k)=
\begin{cases}a^mt^m/m! &\text{if $b=0$, $m_1\geq 0, \dots, m_k\geq 0$,}\\
b^{m}t^{m}/m! &\text{if $a=0$, $m_1\leq 0,\dots, m_k\leq 0$,}\\
0 &\text{otherwise}.
\end{cases}
\end{equation}
\end{cor}
\begin{proof}
It follows from Theorem \ref{thmthw}  and definitions.
\end{proof}

\begin{defn}
Denote $P$ the complete subgraph of $R$ with all nonnegative vertices:

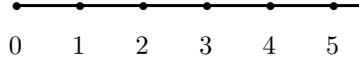
\begin{picture}(410, 25) 
\put(143,18){\circle*{3}}
\put(167,18){\circle*{3}}
\put(191,18){\circle*{3}}
\put(215,18){\circle*{3}}
\put(239,18){\circle*{3}}
\put(263,18){\circle*{3}}
\put(140,0){$0$}
\put(164,0){$1$}
\put(188,0){$2$}
\put(212,0){$3$}
\put(236,0){$4$}
\put(260,0){$5$}
\put(144,18){\line(1,0){22}}
\put(168,18){\line(1,0){22}}
\put(192,18){\line(1,0){22}}
\put(216,18){\line(1,0){22}}
\put(240,18){\line(1,0){22}}
\put(264,18){\line(1,0){10}}
\end{picture}
\end{defn}

\begin{theorem}\label{thm5}
\begin{equation}\label{27}
c_n(P,k \rightarrow l)=\begin{cases}\dbinom{n}{\frac{n+k-l}{2}}-\dbinom{n}{\frac{n-k-l}{2}-1} 
&\text{if $n+k-l$ is even,}\\ 0 &\text{otherwise.}\end{cases}
\end{equation}
\end{theorem}
\begin{proof}
Let $k,l\geq 0$.
Each walk from $k$ to $l$ along the edges of $R$ is a walk along the edges of $P$, except 
the walks coming through $-1$. It means that 
\begin{equation}\label{28}
c_n(P, k \rightarrow l)=c_n(R, k \rightarrow l)-c_n(R,k \rightarrow (-1) \rightarrow l)
=c_n(R,l-k)-c_n(R,-1-k,l+1) . 
\end{equation}
Using Theorems \ref{thm3}  and \ref{thmth}, we get (\ref{27}).
\end{proof}

\begin{cor}\label{corp1}
\begin{equation}\label{eqcorp1}
c(P,k\rightarrow l)=I_{|l-k|}(2t)-I_{l+k+2}(2t) .
\end{equation}
\end{cor}
\begin{proof}
It follows from (\ref{28})  and Corollary \ref{corth}.
\end{proof}

\begin{cor}\label{corp2}
\begin{equation}\label{eqcorp2}
c(P,0\rightarrow m)=\frac{t^m}{m!}\thinspace {_0F_1}(;m+2;t^2)=(m+1)I_{m+1}(2t)/t .
\end{equation}
\end{cor}
\begin{proof}
It follows from Corollary \ref{corp1}  and recurrent relations between Bessel functions.
Or from Theorem \ref{thm5}  and the equality 
\begin{equation}\label{eqbin}
\dbinom{n}{\frac{n-m}{2}}-\dbinom{n}{\frac{n-m}{2}-1}=\frac{m+1}{n+1}\dbinom{n+1}
{\frac{n-m}{2}} .
\end{equation}
\end{proof}

\begin{cor}\label{corp3}
\begin{equation}
\sum_{i=k}^l(n-2i)\binom{n}{i}=n\left( \binom{n-1}{l}-\binom{n-1}{k-1}\right) .
\end{equation}
\end{cor}
\begin{proof}
It follows from changing $n+1$ to $n$ in \eqref{eqbin}, changing $(n-m)/2$ to $i$ and 
summation.
\end{proof}

\begin{defn}\label{Pw}
Denote $P(a,b)$ the complete weighted subgraph of $R(a,b)$ with the same set of vertices as 
$P$.
\end{defn}

\begin{theorem}\label{thm6}
\begin{equation}\label{29}
c_n(P(a,b),k \rightarrow l)=
\begin{cases}\left(  \dbinom{n}{\frac{n+k-l}{2}}-\dbinom{n}{\frac{n-k-l}{2}-1}\right)
a^{(n-k+l)/2}b^{(n+k-l)/2} 
&\text{if $n+k-l$ is even,}\\ 0 &\text{otherwise.}\end{cases}
\end{equation}
\end{theorem}
\begin{proof}
Each walk from $k$ to $l$ 
has the same weight $a^{(n-k+l)/2}b^{(n+k-l)/2}$. Multiplying it on the number of walks stated 
in Theorem \ref{thm5}, we obtain (\ref{29}).
\end{proof} 

\begin{cor}\label{corPw1}For $ab\neq 0$,
\begin{equation}\label{eqcorPw1}
c(P(a,b),k\rightarrow l)=(I_{|l-k|}(2\sqrt{ab}\thinspace t)-I_{l+k+2}(2\sqrt{ab}\thinspace t))
a^{(l-k)/2}b^{(k-l)/2} .
\end{equation}
For $ab=0$,
\begin{equation}\label{eqcorPw12}
c(P(a,b),k\rightarrow l)=\begin{cases}a^{l-k}t^{l-k}/(l-k)! 
&\text{if $b=0$  and $k\leq l$,}\\
b^{k-l}t^{k-l}/(k-l)! 
&\text{if $a=0$  and $k\geq l$,}\\
0 &\text{otherwise.}\end{cases}
\end{equation}
\end{cor}
\begin{proof}
It follows from the weighted variant of (\ref{28})  and Corollary \ref{corth}.
\end{proof}

\begin{cor}\label{corPw2}
For $ab\neq 0$,
\begin{equation}\label{eqcorPw2}
c(P(a,b),0\rightarrow m)=\frac{a^mt^m}{m!}\thinspace {_0F_1}(;m+2;abt^2)=
(m+1)a^{m/2}b^{-m/2}I_{m+1}(2\sqrt{ab}\thinspace t)/t .
\end{equation}
For $ab=0$,
\begin{equation}\label{eqcorPw22}
c(P(a,b),0\rightarrow m)=\begin{cases}a^mt^m/m! 
&\text{if $b=0$,}\\
1 &\text{if $a=0$  and $m=0$,}\\
0 &\text{otherwise.}\end{cases}
\end{equation}
\end{cor}
\begin{proof}
It follows from Corollaries \ref{corp1}, \ref{corp2}  and \ref{corPw1}. 
\end{proof}

\section{Products of linear graphs}\label{sec3} 

\begin{defn}\label{defr^2}
Denote $R^2=R\times R$ the Cartesian square of $R$:

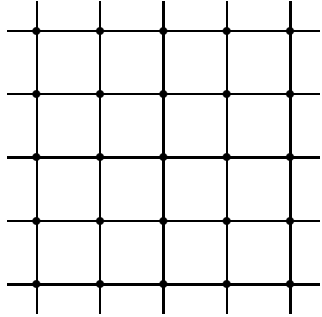
\begin{picture}(410, 125)
\put(155,11){\circle*{3}}
\put(155,10){\line(0,-1){10}}
\put(179,11){\circle*{3}}
\put(179,10){\line(0,-1){10}}
\put(203,11){\circle*{3}}
\put(203,10){\line(0,-1){10}}
\put(227,11){\circle*{3}}
\put(227,10){\line(0,-1){10}}
\put(251,11){\circle*{3}}
\put(251,10){\line(0,-1){10}}
\put(144,11){\line(1,0){10}}
\put(156,11){\line(1,0){22}}
\put(180,11){\line(1,0){22}}
\put(204,11){\line(1,0){22}}
\put(228,11){\line(1,0){22}}
\put(252,11){\line(1,0){10}}

\put(155,35){\circle*{3}}
\put(155,34){\line(0,-1){22}}
\put(179,35){\circle*{3}}
\put(179,34){\line(0,-1){22}}
\put(203,35){\circle*{3}}
\put(203,34){\line(0,-1){22}}
\put(227,35){\circle*{3}}
\put(227,34){\line(0,-1){22}}
\put(251,35){\circle*{3}}
\put(251,34){\line(0,-1){22}}
\put(144,35){\line(1,0){10}}
\put(156,35){\line(1,0){22}}
\put(180,35){\line(1,0){22}}
\put(204,35){\line(1,0){22}}
\put(228,35){\line(1,0){22}}
\put(252,35){\line(1,0){10}}

\put(155,59){\circle*{3}}
\put(155,58){\line(0,-1){22}}
\put(179,59){\circle*{3}}
\put(179,58){\line(0,-1){22}}
\put(203,59){\circle*{3}}
\put(203,58){\line(0,-1){22}}
\put(227,59){\circle*{3}}
\put(227,58){\line(0,-1){22}}
\put(251,59){\circle*{3}}
\put(251,58){\line(0,-1){22}}
\put(144,59){\line(1,0){10}}
\put(156,59){\line(1,0){22}}
\put(180,59){\line(1,0){22}}
\put(204,59){\line(1,0){22}}
\put(228,59){\line(1,0){22}}
\put(252,59){\line(1,0){10}}

\put(155,83){\circle*{3}}
\put(155,82){\line(0,-1){22}}
\put(179,83){\circle*{3}}
\put(179,82){\line(0,-1){22}}
\put(203,83){\circle*{3}}
\put(203,82){\line(0,-1){22}}
\put(227,83){\circle*{3}}
\put(227,82){\line(0,-1){22}}
\put(251,83){\circle*{3}}
\put(251,82){\line(0,-1){22}}
\put(144,83){\line(1,0){10}}
\put(156,83){\line(1,0){22}}
\put(180,83){\line(1,0){22}}
\put(204,83){\line(1,0){22}}
\put(228,83){\line(1,0){22}}
\put(252,83){\line(1,0){10}}

\put(155,107){\circle*{3}}
\put(155,106){\line(0,-1){22}}
\put(155,108){\line(0,1){10}}
\put(179,107){\circle*{3}}
\put(179,106){\line(0,-1){22}}
\put(179,108){\line(0,1){10}}
\put(203,107){\circle*{3}}
\put(203,106){\line(0,-1){22}}
\put(203,108){\line(0,1){10}}
\put(227,107){\circle*{3}}
\put(227,106){\line(0,-1){22}}
\put(227,108){\line(0,1){10}}
\put(251,107){\circle*{3}}
\put(251,106){\line(0,-1){22}}
\put(251,108){\line(0,1){10}}
\put(144,107){\line(1,0){10}}
\put(156,107){\line(1,0){22}}
\put(180,107){\line(1,0){22}}
\put(204,107){\line(1,0){22}}
\put(228,107){\line(1,0){22}}
\put(252,107){\line(1,0){10}}
\end{picture}
\end{defn}

\begin{lemma}\label{npl}
$R^2\simeq R\times_0 R\simeq R\times_1 R$.
\end{lemma}
\begin{proof}
To make this evident, it is enough to rotate the picture above at $45^\circ$ 
counterclockwise. 
\end{proof}

\begin{theorem}\label{thm7}
\begin{equation}\label{30}
c_n(R^2, (m_1,m_2))=\begin{cases}
\dbinom{n}{\frac{n-m_1-m_2}{2}} \dbinom{n}{\frac{n+m_1-m_2}{2}} 
&\text{if $n-m_1-m_2$ is even,}\\ 0 &\text{otherwise.}
\end{cases}
\end{equation}
\end{theorem}
\begin{proof}
It follows from Lemma \ref{npl}, Theorem \ref{thbi1}  and Theorem 
\ref{thm3}. I'll give a direct proof as well.
Consider two projections:
\begin{equation}\label{31}
\begin{split}
&\Pr\!{_1}: R^2\rightarrow R,\quad (m_1,m_2)\mapsto (m_1+m_2) ,\\
&\Pr\!{_2}: R^2\rightarrow R,\quad (m_1,m_2)\mapsto (m_2-m_1) .
\end{split}
\end{equation}
They project each $R^2$-walk from $0$ to $(m_1, m_2)$ to two $R$-walks from $0$ to 
$(m_1+m_2)$  and from $0$ to $(m_2-m_1)$, correspondingly. Note that for each vertex 
of $R^2$, both the projections have the same parity, because 
\begin{equation}\label{31a}
m_1+m_2\equiv m_2-m_1 \mod 2 .
\end{equation} 
Conversely, each pair 
of vertices $u, v$ of $R$, having the same parity, 
uniquely determines the vertex $(m_1, m_2)$ of $R^2$, such that 
\begin{equation}\label{32}
\Pr\!{_1}(m_1, m_2)=u,\quad \Pr\!{_2}(m_1, m_2)=v .
\end{equation}
In fact, 
\begin{equation}\label{33}
m_1=\frac{u-v}{2},\quad m_2=\frac{u+v}{2} .
\end{equation}
For 
\begin{equation}\label{34}
u'=u\pm 1,\quad v'=v\pm 1,
\end{equation}
where the signs can be the same, or different, we have 
\begin{equation}\label{35}
m'_1=m_1,\quad m'_2=m_2\pm 1,
\end{equation}
if the signs were the same, or 
\begin{equation}\label{36}
m'_1=m_1\pm 1,\quad m'_2=m_2,
\end{equation}
if the signs were different. It means that $((m_1,m_2), (m'_1,m'_2))$ is an edge of $R^2$.
In other words, each pair of edges of $R$ having the same 
parity, uniquely determines the edge of $R^2$, the first ( and second) projection of which 
coincides with the first (or second, correspondingly) edge of that pair. Therefore, each 
pair of $R$-walks, from $0$ to $(m_1+m_2)$  and from $0$ to $(m_2-m_1)$, uniquely determines 
a $R^2$-walk from $0$ to $(m_1, m_2)$. We constructed a bijection between the set of 
$R^2$-walks from $0$ to $(m_1, m_2)$ and the direct product of the sets of $R$-walks 
from $0$ to $(m_1+m_2)$  and $R$-walks from $0$ to $(m_2-m_1)$. Now, Theorem \ref{thm3} gives us 
(\ref{30}).
\end{proof}

\begin{cor}\label{cor6}
For $m_1\geq 0, m_2\geq 0$,
\begin{equation}\label{37}
c(R^2, (\pm m_1, \pm m_2))=I_{m_1}(2t)I_{m_2}(2t)=\sum_{i=0}^{\infty}
\binom{m_1+m_2+2i}{m_1+i}\frac{t^{m_1+m_2+2i}}{i!\thinspace (m_1+m_2+i)!}.
\end{equation}
\end{cor}
\begin{proof}
The first equality follows from Theorem \ref{thtran}  and Corollary \ref{cor4}. The equality 
of the first  and the third item of (\ref{37}) follows from Theorem \ref{thm7}.
\end{proof}

\begin{cor}\label{cor7}
For $a\geq 0, b\geq 0$,
\begin{equation}\label{38}
\sum_{i=0}^k \binom{2k+a+b}{i\quad i+a\quad k-i\quad k-i+b} 
=\binom{2k+a+b}{k}\binom{2k+a+b}{k+a} .
\end{equation}
\end{cor}
\begin{proof}After noticing that 
\begin{equation}\label{39}
\binom{2k+a+b}{i\quad i+a\quad k-i\quad k-i+b}=
\binom{2k+a+b}{2i+a}\binom{2i+a}{i}\binom{2k-2i+b}{k-i} ,
\end{equation}
it follows from (\ref{4}), Theorems \ref{thm3}  and \ref{thm7} for $G_1=G_2=R$ and 
\begin{equation} 
u_1=v_1=0,\quad v_1=m_1=a,\quad v_2=m_2=b,\quad n=2k+a+b,\quad j=2i+a .
\end{equation}
There is also another bijective proof of \eqref{38}, without walks. Notice, that the number 
in the right hand side of \eqref{38} is equal to the number of ways of choosing two subsets 
containing $k$  and $k+b$ elements, of the set containing $2k+a+b$ elements. We can choose 
these subsets, choosing first $i$ elements that will be their intersection, then $k-i$ the 
other elements of the first set and $k+b-i$ the other elements of the second set. It gives 
us exactly the sum in the left hand side of \eqref{38}.
\end{proof}

Consider the graph $R\times P$ (the square lattice in the upper-half plane): 

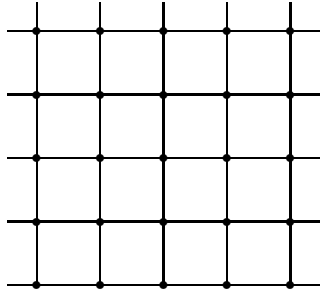
\begin{picture}(410, 116)
\put(155,2){\circle*{3}}
\put(179,2){\circle*{3}}
\put(203,2){\circle*{3}}
\put(227,2){\circle*{3}}
\put(251,2){\circle*{3}}
\put(144,2){\line(1,0){10}}
\put(156,2){\line(1,0){22}}
\put(180,2){\line(1,0){22}}
\put(204,2){\line(1,0){22}}
\put(228,2){\line(1,0){22}}
\put(252,2){\line(1,0){10}}

\put(155,26){\circle*{3}}
\put(155,25){\line(0,-1){22}}
\put(179,26){\circle*{3}}
\put(179,25){\line(0,-1){22}}
\put(203,26){\circle*{3}}
\put(203,25){\line(0,-1){22}}
\put(227,26){\circle*{3}}
\put(227,25){\line(0,-1){22}}
\put(251,26){\circle*{3}}
\put(251,25){\line(0,-1){22}}
\put(144,26){\line(1,0){10}}
\put(156,26){\line(1,0){22}}
\put(180,26){\line(1,0){22}}
\put(204,26){\line(1,0){22}}
\put(228,26){\line(1,0){22}}
\put(252,26){\line(1,0){10}}

\put(155,50){\circle*{3}}
\put(155,49){\line(0,-1){22}}
\put(179,50){\circle*{3}}
\put(179,49){\line(0,-1){22}}
\put(203,50){\circle*{3}}
\put(203,49){\line(0,-1){22}}
\put(227,50){\circle*{3}}
\put(227,49){\line(0,-1){22}}
\put(251,50){\circle*{3}}
\put(251,49){\line(0,-1){22}}
\put(144,50){\line(1,0){10}}
\put(156,50){\line(1,0){22}}
\put(180,50){\line(1,0){22}}
\put(204,50){\line(1,0){22}}
\put(228,50){\line(1,0){22}}
\put(252,50){\line(1,0){10}}

\put(155,74){\circle*{3}}
\put(155,73){\line(0,-1){22}}
\put(179,74){\circle*{3}}
\put(179,73){\line(0,-1){22}}
\put(203,74){\circle*{3}}
\put(203,73){\line(0,-1){22}}
\put(227,74){\circle*{3}}
\put(227,73){\line(0,-1){22}}
\put(251,74){\circle*{3}}
\put(251,73){\line(0,-1){22}}
\put(144,74){\line(1,0){10}}
\put(156,74){\line(1,0){22}}
\put(180,74){\line(1,0){22}}
\put(204,74){\line(1,0){22}}
\put(228,74){\line(1,0){22}}
\put(252,74){\line(1,0){10}}

\put(155,98){\circle*{3}}
\put(155,97){\line(0,-1){22}}
\put(155,99){\line(0,1){10}}
\put(179,98){\circle*{3}}
\put(179,97){\line(0,-1){22}}
\put(179,99){\line(0,1){10}}
\put(203,98){\circle*{3}}
\put(203,97){\line(0,-1){22}}
\put(203,99){\line(0,1){10}}
\put(227,98){\circle*{3}}
\put(227,97){\line(0,-1){22}}
\put(227,99){\line(0,1){10}}
\put(251,98){\circle*{3}}
\put(251,97){\line(0,-1){22}}
\put(251,99){\line(0,1){10}}
\put(144,98){\line(1,0){10}}
\put(156,98){\line(1,0){22}}
\put(180,98){\line(1,0){22}}
\put(204,98){\line(1,0){22}}
\put(228,98){\line(1,0){22}}
\put(252,98){\line(1,0){10}}
\end{picture}

\begin{cor}\label{RP1}
$c(R\times P, (k_1,k_2)\rightarrow (l_1,l_2))=I_{|l_1-k_1|}(2t)
(I_{|l_2-k_2|}(2t)-I_{l_2+k_2+2}(2t))$.
\end{cor}
\begin{proof}
It follows from Theorem \ref{thm1}  and Corollaries \ref{cor4}  and \ref{corp1}.
\end{proof}

\begin{cor}\label{RP2}
$c_n(R\times P, (k_1,k_2)\rightarrow (l_1,l_2))=$
\begin{equation}
\dbinom{n}{\frac{n-(l_1-k_1)-|l_2-k_2|}{2}}\dbinom{n}{\frac{n+(l_1-k_1)-|l_2-k_2|}{2}}-
\dbinom{n}{\frac{n-(l_1-k_1)-(l_2+k_2+2)}{2}}\dbinom{n}{\frac{n+(l_1-k_1)-(l_2+k_2+2)}{2}} ,
\end{equation}
if $n-(l_1-k_1)-(l_2-k_2)$ is even and $0$ otherwise.
\end{cor}
\begin{proof}
It follows from Theorem \ref{thm7}  and Corollaries \ref{cor6}  and \ref{RP1}.
\end{proof} 

\begin{cor}\label{RP3}
$c(R\times P, (k, 0)\rightarrow (l, m))=(m+1)I_{|l-k|}(2t)I_{m+1}(2t)/t$.
\end{cor} 
\begin{proof}
It follows from Theorem \ref{thm1}  and Corollaries \ref{cor4}  and \ref{corp2}.
\end{proof}

\begin{cor}\label{RP4}
\begin{equation}
c_n(R\times P, (k,0)\rightarrow (l,m))= 
\frac{m+1}{n+1}\dbinom{n+1}{\frac{n-(l-k)-m}{2}}\dbinom{n+1}{\frac{n+(l-k)-m}{2}} ,
\end{equation}
if $n-l+k-m$ is even and $0$ otherwise. 
\end{cor}
\begin{proof}
It follows from Theorem \ref{thm7}  and Corollaries \ref{cor6}  and \ref{RP3}.
\end{proof}

\begin{cor}\label{RP5}
\begin{equation}\label{eqRP5}
\binom{n}{i}\binom{n}{j}-\binom{n}{i-1}\binom{n}{j-1}=\frac{n-i-j+1}{n+1}
\binom{n+1}{i}\binom{n+1}{j} .
\end{equation}
\end{cor}
\begin{proof}
Compare Corollaries \ref{RP2}  and \ref{RP4}.
\end{proof}

\begin{cor}\label{RP6}
\begin{equation}\label{eqRP6}
\sum_{i=k}^l(n-m-2i)\binom{n}{i}\binom{n}{m+i}=n\left( \binom{n-1}{l}\binom{n-1}{m+l}-
\binom{n-1}{k-1}\binom{n-1}{m+k-1} \right) .
\end{equation}
\end{cor}
\begin{proof}
Changing $n+1$ to $n$ in \eqref{eqRP5}, changing $j$ to $m+i$ and summing from $i=k$ to 
$l$, we get \eqref{eqRP6}.
\end{proof} 

Analogously, we can obtain a series of similar formulas for $P\times P$  and other 
products of linear graphs. 

\section{Enumeration of plane symplectic wave graphs}

Recall the definition from \cite{MS}. 

\begin{def}\label{ddd1} 
A plane symplectic wave graph is a graph with the vertices $1,2, \dots, m$, 
each connected component of which is a path of length $\geq 1$ 
(i.\ e. it can't be a point), edges of which can be drawn 
in the plane considered as two half-planes, upper  and lower, glued 
along $\RR$, such that the first edge of each connected component, $\{i_1i_2\}$, 
is drawn on the upper half-plane; each edge $\{i_ji_{j+1}\}$ consequent to the    
edge $\{i_{j-1}i_j\}$ drawn on one of the half-planes, is drawn on another one; 
the last edge of the path, 
$\{i_l i_{l+1}\}$, supposed to be drawn on the upper half-plane; 
we suppose also that 
$i_1<i_2<\dots<i_{l+1}$  and edges of our plane symplectic wave graph 
don't intersect.
\end{def}

\begin{theorem}\label{ttt1}
There exist exactly 
\begin{equation}\label{ttt1eq1}
\frac{(2k)!(2k+2)!\cdot 6}{k!(k+1)!(k+2)!(k+3)!}
\end{equation}
plane symplectic wave graphs with $2k$ vertices.
\end{theorem}
\begin{proof}
It was shown in \cite{MS} that there is $\tilde{f}_{2k}^0(2)$ plane symplectic 
wave graphs with $2k$ vertices, where  $\tilde{f}_{2k}^0(2)$ is the number of 
balanced symplectic lattice words in the alphabet $\{1,2,\overline{1},\overline{2}\}$. 
which is equal $c_{2k}(\Lambda^2P, (0,1)\rightarrow (0,1))$. By Theorem \ref{nt5}, 
\begin{equation}\label{ttt1eq2}
c(\Lambda^2P, (0,1)\rightarrow (0,1))=
\begin{vmatrix}c(P,0\rightarrow 0)&c(P,0\rightarrow 1)\\
c(P,1\rightarrow 0)&c(P,1\rightarrow 1)
\end{vmatrix} .
\end{equation}
Plugging the values of $c(P,i\rightarrow j)$ given by Corollaries \ref{corp1}  and 
\ref{corp2}  and using Corollary \ref{cor6}, we obtain formula \eqref{ttt1eq1}.
\end{proof}

\begin{acn}
I would like to thank Fan Chung Graham, Alexandre Kirillov and 
Herbert Wilf for useful discussions
and my Gorgeous and Brilliant Wife, Bette, for her total support and love.
\end{acn}

\end{document}